\documentclass[letterpaper]{article} 
\usepackage[margin=1in]{geometry}
\usepackage{graphicx}	
\graphicspath{{Figures/}}
\usepackage{pdfpages}
 
\usepackage{enumitem}
\usepackage{amssymb}
\usepackage{empheq}
\usepackage{color}
\usepackage{soul,xcolor}
\usepackage{mathabx}
\usepackage{adjustbox}
\usepackage{centernot}
\usepackage{tabularx}
\usepackage{arydshln}
\usepackage{mathtools}
\usepackage{cases}
\usepackage{booktabs}
\usepackage{mathrsfs}  
\usepackage{csquotes}
\usepackage[outdir=./]{epstopdf}
\usepackage{mathtools}
\usepackage{amsmath}
\usepackage{amsfonts}
\usepackage{nicefrac}
\usepackage{natbib}
\usepackage{hyperref}
\hypersetup{
    colorlinks=true,
  citecolor=cyan, urlcolor = red, linkcolor = red}
\usepackage{nccmath}
\usepackage{pst-node}
\usepackage{csquotes}
\usepackage{arydshln}
\usepackage{multirow}
\usepackage{blkarray}
\usepackage[noend]{algpseudocode}
\usepackage{algorithm}
\usepackage{etoolbox}
\usepackage{soul}
\usepackage{scrextend}
\usepackage{accents}

\newcommand{\bb}{\mathbf}
\newcommand{\bbx}{\mathbf{x}}
\newcommand{\bbu}{\mathbf{u}}
\newcommand{\bbA}{\mathbf{A}}
\newcommand{\bbB}{\mathbf{B}}
\newcommand{\bbK}{\mathbf{K}}
\newcommand{\bbT}{\mathbf{T}}
\newcommand{\mf}{\mathfrak}
\newcommand{\mc}{\mathcal}
\newcommand{\bmcK}{\bs{\mathcal{K}}}
\newcommand{\nn}{\nonumber}
\newcommand\thickbar[1]{\accentset{\rule{.4em}{.8pt}}{#1}}
\usepackage{amssymb}
\newtheorem{theorem}{\textbf{Theorem}}
\newtheorem{lemma}{\textbf{Lemma}}
\newtheorem{proposition}{\textbf{Proposition}}
\newtheorem{remark}{\textbf{Remark}}
\newtheorem{assumption}{\textbf{Assumption}}
\newtheorem{definition}{\textbf{Definition}}
\newtheorem{example}{\textbf{Example}}
\newcommand{\tss}{\textsuperscript}
\newcommand{\tsb}{\textsubscript}
\newcommand{\te}{\text}
\newcommand{\bs}{\boldsymbol}
\newcommand{\sbw}{S_{\text{BW}}}
\newcommand{\scn}{S_{\text{CN}}}
\newcommand{\scnc}{S_{\text{CN}_c}}
\newcommand{\scnr}{S_{\text{CN}_r}}
\newcommand{\nnz}{\texttt{nnz}}
\newcommand{\nrow}{{n_{\text{row}}}}
\newcommand{\ncol}{{n_{\text{col}}}}
\newcommand{\ncc}{n_{\text{cc}}}
\newcommand{\ncp}{n_{\text{cp}}}

\newcommand{\noff}{{\bb{n}_{\bb{off}}}}
\newcommand{\noffs}{{n_{\text{off}}}}

\DeclareMathOperator*{\minimize}{minimize}

\allowdisplaybreaks
\newcounter{alphaforalgorithm}
\title{Optimal Co-Designs of Communication and Control in Bandwidth-Constrained Cyber-Physical Systems} 

\author{Nandini Negi\tss{1,2} and Aranya Chakrabortty\tss{1,3}\\
\tss{1}Electrical $\&$ Computer Engineering, North Carolina State University\\
Email: \tss{2}{nnegi@ncsu.edu}, \tss{3}{aranya.chakrabortty@ncsu.edu} }
\date{}
\begin{document}
\maketitle


\begin{abstract}
  We address the problem of sparsity-promoting optimal control of cyber-physical systems (CPSs) in the presence of communication delays. The delays are categorized into two types - namely, an {\it inter-layer} delay for passing state and control information between the physical layer and the cyber layer, and an {\it intra-layer} delay that operates between the computing agents, referred to here as {\it control nodes} (CNs), within the cyber-layer. Our objective is to minimize the closed-loop $\mc{H}_2$-norm of the physical system by co-designing an optimal combination of these two delays and a sparse state-feedback controller while respecting a given bandwidth cost constraint. We propose a two-loop optimization algorithm for this. Based on the alternating directions method of multipliers (ADMM), the inner loop handles the conflicting directions between the decreasing $\mc{H}_2$-norm and the increasing sparsity level of the controller. The outer loop comprises a semidefinite program (SDP)-based relaxation of non-convex inequalities necessary for closed-loop stability. Moreover, for CPSs where the state and control information assigned to the CNs are not private, we derive an additional algorithm that further sparsifies the communication topology by modifying the row and column structures of the obtained controller, resulting in reassigning the communication map between the cyber and physical layers, and determining which physical agent should send its state information to which CN. Proofs for closed-loop stability and optimality are provided for both algorithms, followed by numerical simulations. 
\end{abstract}


 \section{Introduction}

\par \noindent Over the recent years, sparsity-promoting optimal control has emerged as a key tool for enabling economic control of large-scale cyber-physical systems (CPSs) \citep{hespanha,sinopoli} in both continuous-time \citep{mihailo} and discrete-time settings \citep{geromel1989structural}. The fundamental idea is to minimize the number of communication links needed for control without sacrificing the closed-loop performance of the physical system below a specified threshold. Optimization methods such as alternating directions method of multipliers (ADMM) \citep{mihailo,boyd}, proximal Newton method \citep{wytock}, gradient support pursuit (GraSP) \citep{feier}, and rank-constrained convex optimization \citep{motee}, among others, have been successfully used to achieve this trade-off with applications to a wide range of CPSs such as electric power systems, robotics, transportation networks, and multi-agent control. An underlying assumption behind these designs is that the communication of state and control inputs between the different agents is instantaneous. In reality, however, all practical CPSs will encounter communication delays arising from propagation as well as from routing and queuing. How these conventional sparse control designs would perform in the presence of such delays is still an open question. Our recent work in \cite{negi2020sparsity} showed that inclusion of delays is not just a trivial extension of the conventional algorithms for sparse control, but instead demands an entirely new design approach due to various complex stability constraints that are typical to time-delayed systems \cite{delay2,delay1,delay3,liu2012control,hale}.

\par While the work in \cite{negi2020sparsity} addressed this problem for a specific class of CPSs that operate over peer-to-peer communication, in this paper, we present a new design framework for sparsity-promoting optimal control of linear time-invariant (LTI) systems with feedback delay that are defined over a far more generic cyber-physical architecture. Our approach is inspired by recent advancements in cloud computing, fog computing, and software-defined networking (SDN). The physical layer in our CPS consists of the physical plant that needs to be controlled, including its sensors, estimators, actuators, and other physical devices, while the cyber layer is defined over a cloud computing network consisting of multiple spatially distributed {\it virtual} computing agents, referred to here as control nodes (CNs) \citep{xin2011virtual}. The physical layer sensors collectively measure the instantaneous values of the system state and communicate them to designated CNs through a local area network (LAN). Inside the cloud, the CNs then share that state information through an SDN by following the sparsity pattern of the controller. Upon receiving its respective state information, each CN computes a control input using a linear quadratic regulator (LQR) law and sends that information back to a designated actuator in the physical layer, which then actuates that control input. The feedback loop continues like this over time by continuous interactions between the two layers. Unlike the setup in \cite{negi2020sparsity} where a single delay was used, in this case, we have two distinct delays: (1) \textit{inter-layer} delay $\tau_d$ that arises in the LAN connecting the sensors or actuators in the physical layer to the corresponding CNs in the cyber layer, and (2) \textit{intra-layer} delay $\tau_c$ that arises in the SDN links connecting the CNs across the cyber-layer. Both delays are a function of the respective LAN and SDN bandwidths and the distances over which the corresponding communication links are operating \citep{ugtext}. Given this premise, our primary objectives and contributions are as follows.

\par \textbf{[1]} We first present a new sparse optimal control design that minimizes the closed-loop $\mc{H}_2$ norm of the physical system while at the same time designing the optimal values of $\tau_d$ and $\tau_c$ to reduce the bandwidth cost. We co-optimize the controller and these two delays that are all coupled to each other through complex implicit relationships arising from stability, $\mc{H}_2$ performance, and bandwidth constraints. To handle these dependencies, we develop an algorithm (Algorithm \ref{MainAlgorithm}) with two hierarchical loops. The outer loop designs the two delays and finds a corresponding stabilizing controller by sequentially relaxing the non-linear matrix equations required for the co-design. The inner loop sparsifies this controller while minimizing the closed-loop $\mc{H}_2$-norm. Our results show that the relative magnitudes of $\tau_c$ and $\tau_d$ for achieving the optimal $\mc{H}_2$-norm can be notably different depending on the plant dynamics. 
    \smallskip
\par \textbf{[2]} For the case when preserving privacy of the information handled by the CNs is not an issue, we provide a strategy for reassignment of the states and the control inputs by manipulating the block-wise row and column structures of the sparse controller obtained from Algorithm \ref{MainAlgorithm}. This reassignment changes the two delays, resulting in a subsequent change in the $\mc{H}_2$ performance. We propose a series of algorithms, collectively referred to as Algorithm \hyperlink{algo:2}{2}, which further minimize the closed-loop $\mc{H}_2$ norm under the constraint that the computation overhead of the CNs remains below that for Algorithm \ref{MainAlgorithm}. We derive the conditions under which such an optimal reassignment exists.

\par Note that our problem is fundamentally different from the conventional bandwidth allocation and delay assignment problems reported in the literature of computer networking \cite{net1}, \cite{net2}. The utility functions in these papers are static, and do not include any plant dynamics like ours. { We illustrate the effectiveness of our algorithms using simulations in Sec. \ref{sec:simulations1} and \ref{sec:simulations2}. These simulations highlight the impacts of delays and sparsity on $\mc{H}_2$-performance and provide important insights on their interdependencies. 

\par Some preliminary results on this topic have been presented in our recent conference paper \cite{negiacc2021}, where we used a simplified relationship between delay and sparsity to satisfy only the bandwidth cost while minimizing the closed-loop $\mc{H}_2$ norm. The results in this paper, however, are significantly extended, in comparison. We use a more practical delay versus sparsity relationship in our problem formulation, and minimize both the bandwidth and the cost of computation overhead for the CNs. Moreover, the concept of CN reassignment and the related algorithms developed in Sec. \ref{sec:topologydesign} are also added as entirely new contributions.} The rest of the paper is organized as follows. Sec. \ref{sec:problemformulation} states the problem formulation followed by Sec. \ref{sec:proposedcodesign} that describes the proposed co-design of the delays. Sec. \ref{sec:problemsetupinadmm} introduces the two-loop algorithm to solve the problem, followed by the corresponding simulation examples in Sec. \ref{sec:simulations1}. Sec. \ref{sec:topologydesign} introduces the reassignment problem via topology design and derives the corresponding algorithms, followed by simulations in Sec. \ref{sec:simulations2} and conclusion in Sec. \ref{sec:conclusion}. Finally, the proofs of all lemmas, theorems, and propositions are listed in the Appendix.
\par \textbf{Notations:} $\mathbb{R}$, $\mathbb{Z}$ and $\mathbb{N}_n$ are the set of real numbers, integers and natural numbers from $1$ to $n$. $U(0,1)$ is the continuous uniform distribution over [0,1]. The natural order of $i$ refers to ascending order of the indices $i$. $\bbA^T$, $\te{Tr}(\bbA)$ and $\lambda_{max}(\bbA)$ represent the transpose, the trace and the maximum eigenvalue of $\bbA$. $\bbA \otimes \bbB$ and $\bbA\circ \bbB$ represent Kronecker and Hadamard product between $\bbA$ and $\bbB$ respectively. $\bbA^{'}(\bbB)$  represents differentiability of $\bbA$ w.r.t $\bbB$. A permutation matrix is obtained from permuting the rows and columns of $n\times n$ identity matrix $\bb{I}_n$. The rowgroups and colgroups of an $\mc{N}\times \mc{N}$ block matrix $\bbA\in\mathbb{R}^{n\times m}$ refers to partitioning $n$ and $m$ into separate collection of $\mc{N}$ sets. $\texttt{floor}(x)$ rounds $x$ to the nearest integer $\leq x$. The notation $\bb{B}=\texttt{Reshape}(\bbA,[p,q])$ is used to reshape an $\bbA\in\mathbb{R}^{m\times n}$ in row-traversing order to another matrix $\bbB\in\mathbb{R}^{p\times q}$, provided $pq=mn$. 

\begin{figure}[t]
\centering
\includegraphics[scale=0.6]{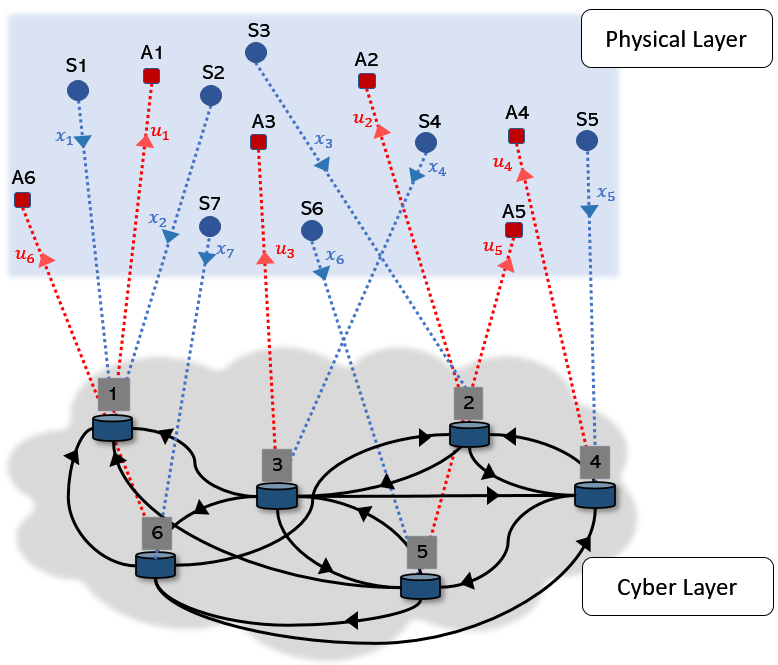}\\
\includegraphics[scale=0.7]{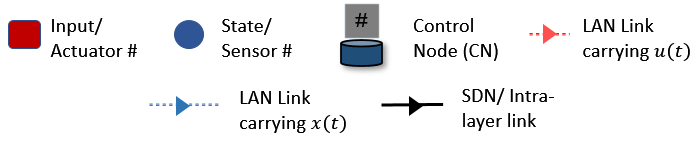}
\caption{Sample CPS schematic showing physical and cyber layers with the associated delays.}
\label{fig:system}
\end{figure}

\section{Problem Formulation}
\label{sec:problemformulation}
\subsection{State Feedback with Communication Delays}
\label{subsec:statefeedbackwithcommdelays}
\noindent Consider a LTI system with the following dynamics:
\begin{equation}
\dot{\bb{x}}(t)=\bbA \bb{x} (t)+ \bb{B} \bb{u}(t) + \bb{B}_w \bb{w}(t), \label{delay-free}
\end{equation}
where $\bb{x} \in \mathbb{R}^{n}$ is the state, $\bb{u} \in \mathbb{R}^{m}$ is the control, and $\bb{w} \in \mathbb{R}^r$ is the exogenous input, with the corresponding matrices $\bbA \in \mathbb{R}^{n\times n}$, $\bb{B}\in\mathbb{R}^{n\times m}$, and $\bb{B}_w\in \mathbb{R}^{n\times r}$. We design a state-feedback controller, ideally represented as $\bb{u}(t)=-\bb{K}\bb{x}(t)$. However, due to limited bandwidth availability, the controller includes finite delays in the feedback. The architecture of the closed-loop system consists of state information $\bb{x}(t)$ being sent from the sensors in the physical layer to the $\mc{N}$ CNs located in a virtual cloud, the CNs sharing this state information with each other and computing the control input $\bbu(t)$ in a distributed way, and finally these control signals being transmitted back to the actuators in the physical layer. The exact CPS model to carry out these three executions is described as follows.
\begin{enumerate}[label=\textbf{A.\arabic*}]
\item \label{A1} Every CN $i$ is associated with $n_i\geq 1$ unique scalar states with $\sum_{i=1}^{\mc{N}} n_i = n$. The subscripts of these states are represented by the set $\mathfrak{x}_i\in\mathbb{Z}^{n_i}$ with $\bigcup_{i=1}^\mc{N} \mf{x}_i = \mathbb{N}_n$. For example, if CN $i$ is associated with $\bbx_1$ and $\bbx_3$, then $\mathfrak{x}_i=\{1,3\}$. These $n_i$ states are received from the physical-layer through LAN or inter-layer communication links with delay $\frac{\tau_d}{2}>0$.
\item \label{A2} Inside the cloud, each CN $i$ shares its corresponding state $\bbx_l(t-\frac{\tau_d}{2})$, $l\in\bs{\mf{x}}_i$ with the others over point-to-point SDN or intra-layer links with delay $\tau_c>0$. 
\item \label{A3} Each CN $j$ is also associated with $m_j\geq 1$ unique scalar control inputs with $\sum_{j=1}^{\mc{N}} m_j = m$. The subscripts of these inputs are represented by the set $\mathfrak{u}_j\in\mathbb{Z}^{m_j}$ with $\bigcup_{j=1}^\mc{N} \mf{u}_j = \mathbb{N}_m$. For e.g., if CN $j$ is associated with $\bbu_2$ and $\bbu_5$, then $\mathfrak{u}_j = \{2,5\}$. These $m_j$ inputs are calculated by CN $j$ at each time $t$.
\item \label{A4}  Each CN $j$ calculates $\bbu_k(t)$ $=$ $-\sum_{l=1}^n\bbK_{kl}$ $\bbx_l(t-\tau_c - \frac{\tau_d}{2})$ $\forall$ $k\in\bs{\mf{u}}_j$, which are then transmitted back to the physical layer with delay $\frac{\tau_d}{2}$. The total round trip delay is, therefore, $\tau_o=\tau_c + \tau_d$.
\end{enumerate}
\par  A sample CPS with $n=7$ states, $m=6$ control inputs, and $\mc{N}=6$ CNs is shown in Fig. \ref{fig:system}. The control input $\bb{u}(t)$ for the CPS described by \ref{A1}-\ref{A4} can be expressed as:
\begin{align}
& \bb{u}(t)  = - \underbrace{(\bb{K}\circ \bs{\mc{I}}_d)}_{\bb{K}_d} \bb{x}(t-\tau_d) - \underbrace{(\bb{K}\circ \bs{\mc{I}}_o)}_{\bb{K}_o} \bb{x}(t-\tau_o), \label{Ueq}
\end{align}
where $\bs{\mc{I}}_d$ , $\bs{\mc{I}}_o\in\mathbb{R}^{m\times n}$ are binary matrices such that
\begin{align}
  \bs{\mc{I}}_d(i,j) =\begin{cases}
    1, \ \text{If} \ \exists \ q\in\{1,\ldots,\mc{N}\} : i \in \bs{\mf{u}}_q, \   j \in \bs{\mf{x}}_q, \ \\
    0, \ \text{otherwise}
    \end{cases},
\end{align}
and $\bs{\mc{I}}_o$ is the complement of $\bs{\mc{I}}_d$. For instance, $\bs{\mc{I}}_d$ and $\bs{\mc{I}}_o$ for the CPS of Fig. \ref{fig:system} are given as:
\begin{equation}
 \bs{\mc{I}}_d= \begin{bmatrix}
     1 &   &   &   &   &   &  \\
   & 1 & 1 &   &   &   &  \\
   &   &   & 1 &   &   &  \\
   &   &   &   & 1 &   &  \\
   &   &   &   &   & 1 &  \\
   &   &   &   &   &   & 1 
    \end{bmatrix}, \ \bs{\mc{I}}_o = \begin{bmatrix}
  & 1 & 1  & 1  & 1   & 1  &1  \\
  1&  &  & 1  & 1  & 1  & 1 \\
  1 &  1 & 1  &  & 1  &  1 & 1 \\
  1 & 1  & 1  & 1  &  & 1  &1  \\
  1 & 1  &  1 &  1 &  1 &  & 1 \\
  1 &  1 & 1  & 1  & 1  & 1  &  \\
    \end{bmatrix}
\end{equation}
The closed-loop system of \eqref{delay-free}-\eqref{Ueq} can be written as:
\begin{align}
&\dot{\bb{x}}(t) = \bbA \bb{x}(t) - \bb{B} \textbf{K}_d \bb{x}(t-\tau_d) - \bb{B} \bb{K}_o \bb{x}(t-\tau_o) + \bb{B}_w \bb{w}(t), \nn \\
&\bb{z}(t) = \bb{C} \bb{x}(t) + \bb{D} \bb{u}(t), \ \bb{C} =[\bb{Q}^{\nicefrac{1}{2}}, \bs{0} ]^T, \ \bb{D}=[\bs{0}, \bb{R}^{\nicefrac{1}{2}}]^T,  \label{delayed}
\end{align}
where $\bb{z}(t)$ is the measurable output, $\bb{Q} \succeq 0$ and $ \bb{R}\succ 0$. We make the standard assumption that $(\bbA,\bb{B})$ and $(\bbA,\bb{Q}^{\nicefrac{1}{2}})$ are stabilizable and detectable, respectively \cite[Sec. II]{mihailo}. Before proceeding, we introduce the following three terms that will be used frequently over the rest of the paper.
\begin{definition}
\label{def:topology}
{\it     
Let us define tuples $\bs{\mf{X}}:=(\bs{\mf{x}}_i)$ and $\bs{\mf{U}}:=(\bs{\mf{u}}_i)$ that respectively represent the state and input indices arranged in the natural order of $i$. For e.g., if $\bs{\mf{x}}_1=\{3,4\}$ and $\bs{\mf{x}}_2=\{2,1\}$, then $\bs{\mf{X}} = (3,4,2,1)$. These tuples provide the order in which states and inputs are allocated to the CNs; for instance, the first $n_1$ ($m_1$) values in $\bs{\mf{X}}$ ($\bs{\mf{U}}$) correspond to the states (inputs) associated with the first CN, the next $n_2$ values with the second CN, and so on. The topology of a CPS, i.e., the state and control inputs associated with all the CNs in the cyber-layer, is defined for $\mc{N}\in [2,\mc{N}_m=\min(m,n)]$\footnote{Since each CN must be associated with at least $1$ state and $1$ control input as given in \ref{A1} and \ref{A3}, the maximum number of CNs in the cyber-layer can be at most $\min(m,n)$.} CNs by the tuple $\bb{T}:=(\bs{\mf{X}},\bs{\mf{U}},\bs{\mf{n}},\bs{\mf{m}})$, where $\bs{\mf{n}} := \{n_i\}$ and $\bs{\mf{m}} = \{m_i\}$.}
\end{definition}
\begin{definition}
\label{def:propdelay}
{\it The propagation delay $\tau_{cpr}$ ($\tau_{{dpr}}$) in the intra-layer (inter-layer) link is the delay arising due to the physical distance between the source and the destination. }
\end{definition}
\begin{definition}
\label{def:transdelay}
{\it The transmission delay $\tau_{ctr}$ ($\tau_{dtr}$) in any intra-layer (inter-layer) link, discussed shortly in Sec. \ref{subsec:problemsetup}, is the delay arising from routing and queuing. This delay is a function of the SDN (LAN) bandwidth, and the total number of corresponding links \cite{ugtext}. In our CPS setting, since both the bandwidth and the number of links can be selected a priori, we will use $\tau_{ctr}$ and $\tau_{dtr}$ as design variables. }
\end{definition}

\subsection{Problem Setup}
\label{subsec:problemsetup}
\par \noindent Our goal is to design a $\bbK$ that minimizes the $\mc{H}_2$-norm of the transfer function from $\bb{w}(t)$ to $\bb{z}(t)$ for the time-delayed LTI system (\ref{delayed}). In general, 
the $\mc{H}_2$-performance of \eqref{delayed} will be worse than that of the delay-free system \cite[Sec. 5.6]{stability}. Therefore, reducing both $\tau_d$ and $\tau_c$ will improve the $\mc{H}_2$-performance. Of course, the trivial solution would be to use $\tau_d=\tau_c=0$, which is not practical as that would require infinite bandwidth and all the link lengths to be zero (Def. \ref{def:propdelay}, \ref{def:transdelay}). We next define constraints on bandwidth and CN computation overhead costs that lower bound the two delays $\tau_d$ and $\tau_c$.
\par -- \textbf{CN cost $\scn$} is the sum of cost of renting and computation overhead of each CN (i.e., number of states and control inputs handled by the CN). This cost is constant for a given topology. In Sec. \ref{sec:topologydesign}, we use topology as a design variable, and accordingly, $\scn$ varies. We provide a mathematical definition of $\scn$ when we come to Sec. \ref{sec:topologydesign}.
\par -- \textbf{Bandwidth cost $\sbw$} is the cost of allocating the total bandwidth to the LAN (inter-layer) and SDN (intra-layer). Let the combined bandwidth of the LAN and SDN links be denoted by $b_{cp}$ and $b_{cc}$, respectively. Then, $\sbw$ can be written as:
\begin{align}
 \sbw =  m_{cp} b_{cp} + m_{cc} b_{cc}, \label{maincost}
\end{align} 
where $m_{cp}$ and $m_{cc}$ are the respective dollar costs for renting LAN and SDN links.
\smallskip
\par The total cost $S$ is $\scn + \sbw$. The bandwidths $b_{cp}$ and $b_{cc}$ are divided according to
their respective number of links as follows.
    \subsubsection{Division of \texorpdfstring{$b_{cc}$}{bcc}}
    \label{subsubsec:sdn} 
        For ease of exposition, let us denote $\bmcK(\bbK,\bb{T})\in\mathbb{R}^{m\times n}$ as the block matrix obtained by first permuting the $m$ rows and $n$ columns of $\bbK$ to follow the ordering of $\bs{\mathfrak{U}}$ and $\bs{\mathfrak{X}}$ (Def. \ref{def:topology}), and then partitioning into $\bs{\mf{n}}$ \textit{rowgroups} and $\bs{\mf{m}}$ \textit{colgroups}. This is shown in the following example.
      \begin{example}
    \label{eg:1}
    {\it For the CPS in Fig. \ref{fig:system}, $\bb{T}$ and $\bbK$ are given as:
    \begin{align}
 \bb{T}\equiv  \Big( \bs{\mf{X}} = (1,2,3,4,5,6,7), \ &\bs{\mf{U}} = ( 1,2,3,4,5,6), \      \bs{\mf{n}} = [1,2,1,1,1,1]^T, \    \bs{\mf{m}} = [1,1,1,1,1,1]^T \Big) \\
 & \bbK=\begin{bmatrix}
 a &   &   & b &   & c & d \\
   & e & f &   & g &   & h\\
  & i & j &   & k & \\
  & l &m  &  &  & n \\
&  & o & p & q &  \\
&  & r &   & s & 
 \end{bmatrix}.
    \end{align}
We obtain the block matrix $\bmcK(\bbK,\bb{T})$ as:
  \begin{equation}
   \setlength\aboverulesep{0pt}\setlength\belowrulesep{0pt}
    \setlength\cmidrulewidth{0.5pt}
  \begin{blockarray}{cccccccc}
  {{}_{u}\mkern-1mu\setminus\mkern-1mu{}^{x}}    & {\overbrace{1}^{\bs{\mf{x}}_1}} & \BAmulticolumn{2}{c}{\overbrace{2 \hspace{0.48cm} 3}^{\bs{\mf{x}}_2}} & {\overbrace{4}^{\bs{\mf{x}}_3}}& {\overbrace{5}^{\bs{\mf{x}}_4}}&
 {\overbrace{6}^{\bs{\mf{x}}_5}}&
{\overbrace{7}^{\bs{\mf{x}}_6}}\\ 
\vspace{-0.4cm} \\
 \begin{block}{c[c|cc|c|c|c|c]}
\bs{\mf{u}}_1 \  \begin{cases} 1 \end{cases}  & a & &   &b  &  &c  & d\\
\cmidrule(lr){2-8}
\bs{\mf{u}}_2   \  \begin{cases} 2 \end{cases} &  & e & f &  & g& & h\\
\cmidrule(lr){2-8}
\bs{\mf{u}}_3  \  \begin{cases} 3 \end{cases} &  &  & i & j &  &  k & \\
\cmidrule(lr){2-8}
\bs{\mf{u}}_4   \  \begin{cases} 4 \end{cases} &  &  & l  & m  &  &  & n \\
\cmidrule(lr){2-8}
\bs{\mf{u}}_5  \  \begin{cases} 5 \end{cases} &  &  &  & o & p  & q & \\
\cmidrule(lr){2-8}
\bs{\mf{u}}_6  \  \begin{cases} 6 \end{cases} &  &  &  & r &  & s  & \\
    \end{block}
\end{blockarray}. \label{Keg1}
  \end{equation} }
 \end{example}
\par Given a $\bbK$ for a fixed $\bbT$, an intra-layer link from CN $i$ to $j$ is not needed if calculation of $\bbu_k$, $k\in\bs{\mf{u}}_j$ does not require $\bbx_l$, $l\in\bs{\mf{x}}_i$. This happens if $\bbK_{kl}=0$ $\forall$ $k\in\bs{\mathfrak{u}}_j$, $l\in\bs{\mathfrak{x}}_i$ (See \ref{A4}), i.e., $\bmcK(\bbK,\bb{T})_{j,i}=\bb{0}$. Here, $\bmcK(\bbK,\bb{T})_{j,i}$ represents the $j,i$-th block of the block matrix $\bmcK$. Therefore, the number of outgoing intra-layer links from CN $i$ are the number of non-zero off-diagonal blocks in the $i$-th block column of $\bmcK(\bbK,\bbT)$, denoted by $n_{\text{off}_i}$. { We define $\noff = \{ n_{\text{off}_i}\}$ as the vector of these entries for all the CNs.}
   \setcounter{example}{0}
  \begin{example}
{\it (Contd.) For the $\bmcK$ in \eqref{Keg1}, $n_{\text{off}_1} = 0$ since all the off-diagonal blocks in the first block-column are zero. This means that none of the CNs (except CN $1$ itself) require the state $\bbx_1$ held by CN $1$ to calculate their respective control inputs. The rest of the vector is obtained as $\bb{n}_\bb{off}=[0,2,4,2,3,3]^T$. Thus, there are a total of $14$ intra-layer links.}
  \end{example}
\par \noindent CN $i$ transmits $n_i$ states in each of the $n_{\text{off}_i}$ outgoing intra-layer links. Thus, the intra-layer bandwidth $b_{cc}$ is divided into the total number of channels in all the links, denoted as $\ncc(\bbK,\bbT)$, resulting in the intra-layer delay
\begin{align}
 \tau_{c} =& \underbrace{\kappa \left(\frac{\ncc(\bbK,\bbT)}{b_{cc}}\right)}_{\tau_{ctr}(\bbK,\bbT)}+ \tau_{cpr} (\bb{T}), \  \label{tauc} 
\end{align}
where $\ncc(\bbK,\bbT)=  \ \bs{\mf{n}}^T(\bbT)\bb{n}_\bb{off}(\bmcK(\bbK,\bb{T}))$ and $\kappa$ is a proportionality constant.
   
\subsubsection{Division of \texorpdfstring{$b_{cp}$}{bcp}}
 \label{subsubsec:lan} 
 The uplink for carrying $\bbu_j$ back to the physical layer is not needed if the $j$-th row of $\bbK$ is entirely $0$. Similarly, if the $i$-th column of $\bbK$ is $0$, then $\bbx_i$ is no longer required for calculating any control input, and the corresponding downlink becomes redundant. Thus, $b_{cp}$ is effectively divided into the number of non-zero rows and columns of $\bbK$ denoted by $\nrow(\bbK)$ and $\ncol(\bbK)$, respectively. The delay in the inter-layer links is, therefore, written as:
    \begin{align}
\tau_d =&\underbrace{2 \kappa \left( \frac{\ncp(\bbK)}{b_{cp}} \right)}_{\tau_{dtr}(\bbK)} + \tau_{dpr}(\bb{T}) , \label{taud}
\end{align}
where $\ncp(\bbK)=  \ \nrow(\bbK) + \ncol(\bbK)$.
\begin{remark}
 { \it Note that the definitions in \eqref{tauc} and \eqref{taud} do not involve any subscript for $\tau_c$ and $\tau_d$ indicating that the delays are assumed to be equal across all links in the SDN as well as in the LAN. While this assumption is made to simplify the design, its practical relevance is as follows. For our practical purposes, propagation delay is of the order of $10^{-4}$ \citep{ugtext} due to which we assume it to be equal for all intra-layer and inter-layer links. For the transmission component of the delay, one can assign the per-link bandwidth in a way that the per-link transmission delay becomes equal for all links.}
  \end{remark}
  \medskip
 \par  Using \eqref{maincost}, \eqref{tauc} and \eqref{taud}, we can write the bandwidth cost constraint as:
\begin{align}
\sbw = & \frac{ 2 m_{cp}\ncp(\bbK)}{\underbrace{\tau_d - \tau_{dpr}}_{\tau_{dtr}}} +  \frac{m_{cc}\ncc(\bbK,\bb{T})}{\underbrace{\underbrace{\tau_o-\tau_d}_{\tau_{c}} - \tau_{cpr}}_{\tau_{ctr}}}\leq S_b,  \label{bw2}
\end{align} 
where $S_b>0$ is a mandatory budget that is imposed to prevent infinite bandwidth. To minimize the cost of renting the links and bandwidth, we wish to reduce the number of both LAN (inter-layer) and SDN (intra-layer) links by promoting sparsity in $\bbK$. Our design objectives, therefore, are listed as:
\par \textbf{P1:} \textit{Given a topology $\bbT$, design $\tau_d$, $\tau_o$ and $\bbK$ such that} 
\par $\bullet$ \textit{$\mc{H}_2$-norm of the closed-loop transfer function of \eqref{delayed} from $\bb{w}(t)$ to $\bb{z}(t)$, denoted as $J$, is minimized.}
\par $\bullet$ \textit{The bandwidth cost $\sbw$ satisfies \eqref{bw2} and budget $S_b$, which is assumed to be large enough for the problem to be feasible.}
\par $\bullet$ \textit{Sparsity of $\bbK$ is promoted.}

\par \noindent Given a budget $S_b$ and a topology $\bbT$, \textbf{P1} can be mathematically stated as:
\begin{subequations}
\label{optim1}
\begin{align}
 \textbf{O1}: \ &\underset{\bbK,\tau_d,\tau_o}{\text{minimize}} \ \ \  J(\bbK,\tau_d,\tau_o) + g(\bbK), \\
 &\text{subject to} \ \ \  \text{$\bbK$ stabilizes \eqref{delayed} for $\tau_o$ and $\tau_d$},\\
 & \hspace{1.75cm} \sbw(\tau_d,\tau_o,\bbK) \leq S_b,
  \end{align}
\end{subequations}
where $\sbw$ is given by \eqref{bw2}, and $g(\bbK)$ is a sparsity-promoting function which will be introduced in Sec. \ref{subsection:innerloop}. The closed-form expression of $J$ is derived next.

\subsection{ \texorpdfstring{$\mc{H}_2$}{H2} norm for the Delayed System}
\noindent The delayed system \eqref{delayed} is infinite dimensional. In order to obtain a linear, finite dimensional LTI approximation of \eqref{delayed}, we use the method of spectral discretization given in \cite{spectral}. Since $\tau_o > \tau_d$ in \eqref{delayed}, following \cite{spectral}, we divide $[-\tau_o,0]$ into a grid of $N$ scaled and shifted Chebyshev extremal points
\begin{equation}
\theta_{k+1} = \frac{\tau_o}{2}\left(\cos\left(\frac{(N-k-1)\pi}{N-1}\right) -1\right), \ k=\{0,\ldots,N-1\}, \label{theta}
\end{equation}
such that $\theta_1=-\tau_o$ and $\theta_N=0$. {The choice of $N$ is guided by \cite[Sec. 4]{spectral}}. Let $\bs{\upsilon}(\theta)=\bb{x}(t+\theta)$ denote the $\theta$-shifted state vector. The extended state  $\bs{\eta}$ and the closed-loop state matrix $\bb{A}_{cl}$ can then be written as:
\begin{subequations}
\label{spectral}
\begin{align}
&\bs{\eta} = [\upsilon^T(\theta_1),  \cdots ,\upsilon^T(\theta_N)= x(t)]^T, \ l_j(\theta)=\prod\limits_{m=1,\ m\neq j}^N \frac{\theta -\theta_m}{\theta_j -\theta_m}, \label{lj}\\
& {\bbA_{cl}}_{ij} = \begin{cases}
 \partial_\theta l_{j}(\theta_i) \bb{I}_n, \ \ \ j=1,\ldots,N, \ i=1,\ldots,N-1\\
 -l_N(-\tau_d) \bbB \bbK_d +  \bbA, \ \ \ \  j=N, \ i=N\\
 -l_1(-\tau_d) \bbB \bbK_d - \bbB  \bbK_o, \ \ \ \  j=1, \ i=N\\
 -l_j(-\tau_d) \bbB \bbK_d , \ \ \ \ \ j=2,\ldots,N-1, \ i=N,
\end{cases} \label{atilde}
\end{align}
\end{subequations}
where $\bbK_d = \bbK\circ \bs{\mc{I}}_d$, $\bbK_o = \bbK \circ \bs{\mc{I}}_o$. 
We can separate $\bbA_{cl}$ into three sub-components:
\begin{align}
&\hspace{1.5cm} \bbA_{cl} = \tilde{\bbA} - \bs{\mc{B}} \bbK_o \bb{N}^T_o - \bs{\mc{B}} \bbK_d \bb{N}^T_d, \label{acl}\\
&\bs{\mc{B}}=\bb{M} \bb{B}, \ \bb{M}=[\bb{0},\ldots,\bb{0},\bb{I}_n]^T, \ \bb{N}_o = [\bb{I}_n,\bb{0},\ldots,\bb{0}]^T,
\end{align}
where the first sub-component $\tilde{\bbA}$ is independent of $\bbK_d$ and $\bbK_o$, the second is only dependent on $\bbK_o$, and the third on $\bbK_d$. The explicit expressions for $\tilde{\bbA}$ and $\bb{N}_d$ in terms of $\tau_d$ and $\tau_o$ will be derived in Sec. \ref{subsec:h2vsdesignvariables}. The linear approximation of the closed-loop system \eqref{delayed} becomes:
\begin{subequations}
\label{extendeddelayed}
\begin{align}
&\dot{\bs{\eta}}(t) = \bbA_{cl} \bs{\eta}(t) + \bs{\mc{B}}_w\bb{w}(t), \\
& \bb{z}(t) = \bs{\mc{C}}\bs{\eta}(t), \  \bs{\mc{C}} = [\bb{M}^T\bb{Q}^{\nicefrac{1}{2}},\ -(\bbK_d \bb{N}^T_d + \bbK_o \bb{N}^T_o)^T\bb{R}^{\nicefrac{1}{2}}  ]^T
\end{align}
\end{subequations}
where  $\bs{\mc{B}}_w =\bb{M}\bb{B}_w$. The algebraic Riccati equations and the closed-loop $\mc{H}_2$-norm $J$ can be written as:
\begin{align}
&\bbA_{cl}^T\bb{P}+\bb{P}\bbA_{cl} = -\bs{\mc{C}}^T\bs{\mc{C}} =-\big( \tilde{\bb{Q}} + \tilde{\bb{C}}^T \bb{R} \tilde{\bb{C}}\big), \label{CC}  \\
&\bbA_{cl}\bb{L}+\bb{L}\bbA_{cl}^T = -\bs{\mc{B}}\bs{\mc{B}}^T ,  \label{LL}\\
&J(\bbK,\tau_d,\tau_o) = \text{Tr} ( \bs{\mc{B}}^T\bb{P}\bs{\mc{B}} ) = \text{Tr} (\bs{\mc{C}}\bb{L}\bs{\mc{C}}^T). \label{J}
\end{align}
where $\tilde{\bb{Q}} =\bb{M}\bb{Q} \bb{M}^T$ and $\tilde{\bb{C}}= \bb\bbK_d \bb{N}^T_d + \bbK_o \bb{N}^T_o$.

\section{ Derivation of the gradient of \texorpdfstring{$\mc{H}_2$}{H2} norm}
\label{sec:proposedcodesign}

\noindent Our goal is to design $(\bbK,\,\tau_d,\,\tau_o)$ to minimize $J$. However, from (\ref{CC})-(\ref{J}),  we see that $J$ is a function of $\tilde{\bbA}$ and $\bb{N}_d$, besides $\bbK$. To compute the gradient of $J$ with respect to $(\bbK,\,\tau_d,\,\tau_o)$, it is essential to express $\tilde{\bbA}$ and $\bb{N}_d$ in terms of these three design variables. We begin this section with these derivations as follows.

\subsection{\texorpdfstring{$\mc{H}_2$}{H2} Performance and Design Variables}
\label{subsec:h2vsdesignvariables}
\noindent Recall that the closed-loop state matrix $\bbA_{cl} = \tilde{\bbA} - \bs{\mc{B}} (\bb\bbK_o \bb{N}^T_o$ $ +\bbK_d \bb{N}^T_d  )$. In the next two lemmas, we express $\bbA_{cl}$ as a function of $\tau_o$, $\bbK$ and the {delay} ratio $c=\nicefrac{\tau_d}{\tau_o}$.
\begin{lemma}
\label{lemmatau}
{\it $\tilde{\bbA}$ is a function of $\tau_o$, and can be written as:
\begin{align}
\tilde{\bbA} = \frac{1}{\tau_o}\bs{\Lambda} + \thickbar{\bb{A}}, \ \thickbar{\bb{A}}=\texttt{Diag}(\bb{0},A),
\end{align}
where $\bs{\Lambda}$ is a constant matrix for constant $N$.} $\blacksquare$
\end{lemma}
\begin{lemma}
\label{lemmand}
{\it $\bb{N}_d$ is a function of the ratio $c=\nicefrac{\tau_d}{\tau_o} \in [0,1]$, and can be written as:
\begin{equation}
\bb{N}_d (c) = \left(\bs{\Gamma} \bs{\nu} (c)\right) \otimes \bb{I}_n, \ \bs{\nu} (c) = [c^{N-1} \ c^{N-2} \ \ldots \ c^2 \ c \ 1]^T,
\label{Nddeclare}
\end{equation}
where $\bs{\Gamma} \in \mathbb{R}^{N\times N}$ is a constant matrix for constant $N$.} $\blacksquare$
\end{lemma}
\smallskip
\par Lemmas \ref{lemmatau} and \ref{lemmand} show that {for fixed $N$}, $J$ for the system in \eqref{extendeddelayed} is a function of $\tau_o$ and $c$. Henceforth, all of our analysis for minimizing $J$ will be carried out using $\tau_o$ and $c$, instead of $\tau_o$ and $\tau_d$. This change of variables is invertible, and therefore, there is no loss of generality.

\subsection{Gradient of \texorpdfstring{$\mc{H}_2$}{H2} norm}
In order to minimize $J$, we next derive the gradient of $J$. We define the set of solutions that guarantee closed-loop stability of \eqref{extendeddelayed} as:
\begin{align}
 \mathbb{K}:= \{(\bbK,\tau_o,c): \te{Re}\big(\lambda_{max}(\bbA_{cl})\big) < 0 \}.   \end{align}
Given this definition, we next derive the gradient of closed-loop $\mc{H}_2$ norm $J$ at $\bbK$, $\tau_o$ and $c$ in the following theorem.
\begin{theorem}
\label{Theoremgrad}
 \textit{$J$ in \eqref{J} is differentiable on $\mathbb{K}$. With $\bb{P}$ and $\bb{L}$ obtained from \eqref{CC} and \eqref{LL}, the gradient of $J$ is evaluated as:}
\begin{align}
J'(\tau_o)  = & -\frac{2}{\tau^2_o}\te{Tr} ( \bs{\Lambda}^T\bb{P} \bb{L}), \\
J'(c) = &  2 \te{Tr}( \bb{N}^{'}_d(c) \bbK^T_d \bb{G} \bb{L}), \label{JtauJc} \\
\nabla J (\bbK) = & 2 ( (\bb{G}\bb{L}\bb{N}_d)\circ \bs{\mc{I}}_d + (\bb{G}\bb{L}\bb{N}_o)\circ \bs{\mc{I}}_o),  \label{JK}
\end{align}
where $\bb{G}= \bb{R} (\bb\bbK_d \bb{N}^T_d+ \bbK_o \bb{N}^T_o) - \bs{\mc{B}}^T\bb{P}$ and $\bb{N}^{'}_d=\left(\bs{\Gamma}  \partial \bs{\nu}(c)\right)\otimes \bb{I}_n$.
\end{theorem}
The negative directions of $ J'(c)$ and $J'(\tau_o) $, as derived in Theorem \ref{Theoremgrad}, always point to the trivial solution $c=0, \ \tau_o=0$ which defeats the purpose of designing $\tau_d$ and $\tau_o$. This is because the partial derivatives in \eqref{JtauJc}-\eqref{JK} are derived with the assumption that $\bbK$, $\tau_o$ and $c$ are independent of each other as $\bbK^{'}(\tau_o)$ and $\bbK^{'}(c)$ cannot be computed directly given the implicit dependence of $\bbK$ on $\tau_o$ and $c$. Therefore, it would be incorrect to co-design $c$, $\tau_o$ and $\bbK$ using just the gradient information. Starting from a stabilizing $(\bbK,\tau_o,c)\in\mathbb{K}$, as soon as we change either $\tau_o$ or $c$, we must update $\bbK$ to ensure stability of \eqref{extendeddelayed}. In other words, $(\bbK,\tau_o)$ and $(\bbK,c)$ must be co-designed separately in sequence while holding $c$ and $\tau_o$ as constant in the respective steps.
\subsection{Co-design of Controller and Delays}
 We next describe how equations in \eqref{CC}-\eqref{LL} can be relaxed for each of the two co-designs.
 \smallskip
\par \noindent $\bullet$ \textit{\textbf{Co-design of $(\bbK,\tau_o)$}}
\smallskip
\begin{theorem}
\label{theoremktau}
\textit{Let $\omega_o= \nicefrac{1}{\tau_o}$. Consider a known tuple $(\bbK^*,\omega^*_o,c^*)\in\mathbb{K}$ satisfying \eqref{CC} with a known $\bb{P}^*$ for closed-loop state matrix $\bbA^{*}_{cl}(\bbK^{*},\omega^{*},c^{*})$. Let ${\omega_o} = \omega^*_o + {\Delta \omega}$, $\bbK=\bbK^* + {\Delta \bbK}$, $\bb{P}=\bb{P}^* + {\Delta\bb{P}}$ and ${\alpha}>0$ be obtained as a solution of the following SDP:
\begin{subequations}
\label{theoremktaueq}
\begin{align}
&\bs{\phi}_0 + {\bs{\phi}_1} + {\bs{\psi}_0} + {\alpha} \bb{I} \succeq 0, \\
&|{\Delta \omega}| \leq \zeta_1, \ \|{\Delta\bb{P}} \| \leq \zeta_2, \label{39bc} \\
& {\alpha} \geq 2\zeta_1 \|\bs{\Lambda}^T {\Delta\bb{P}}\| + 2\zeta_2 \|\bs{\mc{B}}{\Delta \tilde{\bb{C}}}\| + \| \bb{R}^{\nicefrac{1}{2}}{\Delta \tilde{\bb{C}}}\|^2, \label{39d}
\end{align}
\end{subequations}
where $\alpha$, $\Delta \bbK$, $\Delta\bb{P}$ and $\Delta \omega$ are the design variables, ${\bs{\phi}_0} = \bbA_{cl}^{*T}{\bb{P}} + \bb{P}\bbA^*_{cl}$, $\bbK^*_d = \bbK^* \circ \bs{\mc{I}}_d$, ${\Delta \bbK_d} = {\Delta \bbK} \circ \bs{\mc{I}}_d$, $\bbK^*_o =\bbK^* \circ \bs{\mc{I}}_o$, ${\Delta \bbK_o }= {\Delta \bbK} \circ \bs{\mc{I}}_o$, $\tilde{\bb{C}}^* =  (\bbK^*_d \bb{N}^T_d +\bbK^*_o \bb{N}_o^T)$, ${\Delta \tilde{\bb{C}}}=({\Delta \bbK_d} \bb{N}^T_d +{\Delta \bbK_o} \bb{N}^T_o)$, ${\bbA_1 }= -\bs{\mc{B}} ({\Delta \tilde{\bb{C}}}) + {\Delta \omega} \bs{\Lambda}$, ${\bs{\phi}_1 }= {\bbA}^T_{{1}}\bb{P}^* +\bb{P}^* {\bbA_1}$, ${\bs{\psi}_0} = \tilde{\bb{Q}} + \tilde{\bb{C}}^{*T} \bb{R} \tilde{\bb{C}}^* + {\Delta \tilde{\bb{C}}}^T \bb{R} \tilde{\bb{C}}^* + \tilde{\bb{C}}^* \bb{R} {\Delta \tilde{\bb{C}}}$  and,  $\zeta_1, \ \zeta_2$ are chosen constants. Then, $(\bbK,\nicefrac{1}{\omega_o},c^{*})$ is a stabilizing tuple for \eqref{extendeddelayed}.} $\blacksquare$
\end{theorem}
\smallskip
\par \noindent $\bullet$ \textit{\textbf{Co-design of $(\bbK,c)$}}
\smallskip
\par Next, consider the co-design step for $(\bbK,c)$. Recall that $\bbA_{cl}$ is a non-linear function of $c \in [0,1]$ through $\bb{N}_d(c)$ as shown in Lemma \ref{lemmand}, and therefore, the exact expression of $\bb{N}_d(c)$ cannot be used while forming the SDP relaxations. To circumvent this problem, we divide $[0,1]$ into $k_c$ sub-intervals $[{c_1},{c_2}],\ldots,[c_{k_c},c_{k_c + 1}]$ with each sub-interval small enough to allow $\bb{N}_d(c)$ to be approximated as an affine function $\hat{\bb{N}}_d(c)$. Let each sub-interval $[c_i,c_{i+1}]$ have an associated vector of affine coefficients $\bs{\chi}^{(i)} \in \mathbb{R}^{N\times 2}$. The approximated function is written as:
\begin{align}
&\hat{\bb{N}}_d (c) = \left(\bs{\chi}^{(i)} [c, 1]^T\right)\otimes \bb{I}_n, \ c\in [c_i,c_{i+1}], \ i=1,\ldots,k_c.
\end{align}
The coefficients can be computed from a linear curve fitting on \eqref{Nddeclare}. Larger the number of sub-intervals $k_c$, lower is the approximation error $\|\hat{\bb{N}}_d - \bb{N}_d\|$. For our simulations in Sec. \ref{sec:simulations1}, we have used $k_c=10$. We next present the SDP relaxation for the co-design of $(\bbK,c)$.
\smallskip
\begin{theorem}
\label{theoremkc}
{\it Consider a known tuple $(\bbK^{*},\tau^{*}_o$, $c^{*})\in\mathbb{K}$ with $c^*\in[c_i,c_{i+1}]$ for some $i\in \{ 1,\ldots,k_c\}$ satisfying \eqref{LL} with a known $\bb{L}^{*}$ for closed-loop state matrix $\bbA^{*}_{cl}(\bbK^{*},\tau^{*}_o,c^{*})$. Let $c=c^{*} + \Delta c$, $\bbK=\bbK^{*} + \Delta \bbK$, $\bb{L}=\bb{L}^{*}$ $+\Delta \bb{L}$ and $\alpha>0$ be a solution of the following SDP:
\begin{subequations}
\label{theoremkceq} 
\begin{align}
& \bs{\phi}_0 + \bs{\phi}_1 + \bs{\mc{B}}\bs{\mc{B}}^T + \alpha \bb{I} \succeq 0,  \label{beta}\\
&c\in [c_i, c_{i+1}], \ \|\Delta \bb{L}\| \leq \beta, \\
& \alpha \geq 2 \beta \|\bs{\mc{B}} (\Delta \bbK_d \bb{N}^T_d (c^{*}) + \Delta \bbK_o \bb{N}^T_o)\| + 2(\beta\mathfrak{S}\|\bs{\mc{B}}\Delta \bbK_d \| +  \beta \| \bs{\mc{B}}\bb{K}^{*}_d \Delta \bb{N}^T_d \| +  \mathfrak{S} \|\bs{\mc{B}}\Delta \bbK_d\|\|\bb{L}^*\| ), \label{alphaconstraint}
\end{align}
\end{subequations}
where $\alpha, \ \Delta \bbK, \ \Delta\bb{P}$ and $\Delta c$ are the design variables, $\Delta \bb{N}_d $ $ = \hat{\bb{N}}_d (c) - \bb{N}_d (c^*)$, $\bs{\phi}_0 = \bbA_{cl}^{*}\bb{L}+\bb{L}\bbA_{cl}^{* T}$, $\bs{\phi}_1 = \bbA_1 \bb{L}^{*} + \bb{L}^{*} \bbA^T_1$, $\bbA_1 = -\bs{\mc{B}} (\bbK^{*}_d \Delta \bb{N}^T_d(c) + \Delta \bbK_d \bb{N}^T_d (c^{*}) + \Delta \bbK_o \bb{N}^{T}_o)$, $\beta>0$ is a chosen constant, and $\mathfrak{S}\geq \|\bb{N}_d(c)\|$. Then, $(\bbK,\tau^{*}_o,c)$ is a stabilizing tuple for \eqref{extendeddelayed}.} $\blacksquare$
\end{theorem}
\par \noindent Starting from a known stabilizing tuple $(\bbK^*,\tau^*,c^*)$, Theorems \ref{theoremktau} and \ref{theoremkc} enable us to co-design stabilizing pairs $(\bbK,\tau_o)$ and $(\bbK,c)$, respectively. We next integrate the bandwidth cost constraint \eqref{bw2} with the SDPs in \eqref{theoremktaueq} and \eqref{theoremkceq}.
\subsection{Incorporating Bandwidth Constraints}
\label{subsec:changeofvariables}
\par \noindent We impose the bandwidth cost constraint \eqref{bw2} as part of \textbf{P1}, which can be rewritten as:
\begin{align}
&\sbw =\frac{ 2m_{cp}\ncp(\bbK) }{c \tau_o - \tau_{dpr}} + \frac{m_{cc}\ncc(\bbK,\bbT)}{\thickbar{c}\tau_o - \tau_{cpr}}\leq S_{b}, \label{constraintchange}
\end{align}
where $\ncp(\bbK) =\nrow(\bbK) + \ncol(\bbK)$, $\ncc(\bbK) = \bs{\mf{n}}^T(\bbT) \noff(\bbK,\bbT)$ and $\thickbar{c}=1-c$. Recall that $\sbw$ is the total bandwidth cost and $S_b$ is the upper bound imposed on it as stated in \textbf{O1}. When \eqref{constraintchange} is imposed on SDPs \eqref{theoremktaueq} and \eqref{theoremkceq}, we obtain an alternative form of \eqref{constraintchange}, which is stated in the next proposition.
\begin{proposition}
\label{propositionkcktau}
{\it  Given the topology $\bb{T}$, with corresponding constant propagation delays $\tau_{dpr}$ and $\tau_{cpr}$, consider a known tuple $(\bbK^{*},\tau^*_o,c^*)\in\mathbb{K}$ with an associated bandwidth cost $\sbw^*\leq S_b$. Denoting $n^*_{cp}$ $=\nrow(\bbK^*)+\ncol(\bbK^*)$ and $n^*_{cc}=\bs{\mf{n}}^T\noff(\bbK^*)$
\footnote{\label{footnote:constanttopology}Since the topology $\bb{T}$ is assumed to be constant, with a slight abuse of notation we write $\noff$ as a function of $\bbK^*$.}, the following statements are true.
\par 1) Keeping $\tau_o=\tau^*_o$, let $c^*$ be perturbed to $c\in$ $\left({\tau_{dpr}}/{\tau_o^*},1-\right.$ $\left.({ \tau_{cpr}}/{\tau_o^*}) \right)$ resulting in a cost $\sbw(c)$. Then, $\delta \sbw(c): = \sbw-\sbw^*\leq 0$ is a convex constraint and written as:
\begin{align}
 &\delta \sbw(c) = (p_1-p_2) c^2 + (q_1-q_2) c + (r_1-r_2) \leq 0, \label{constraintchange2}  
 \end{align}
 where $p_1 = \sbw^* {\tau_o^*}^2$, $p_2 = -\tau_o^*$, $q_1 =  \tau_o^* \left(-\sbw^* (\tau_o^* + \right.$ $\left. ( \tau_{dpr} - \tau_{cpr}) )- 2 m_{cp} n_{cp}^* + m_{cc} 
 \ncc^* \right)$, $q_2 = \tau_o^* \left(\tau_o^* + (  \right.$ $\left. \tau_{dpr}  - \tau_{cpr})  \right)$, $r_1 = (\sbw^*\tau_{dpr} + 2m_{cp} n_{cp}^* )(\tau_o^* - \tau_{cpr}) - m_{cc}\ncc^* \tau_{dpr}$ and $r_2  = -\tau_{dpr} (\tau_o^* - \tau_{cpr})$. The constraint $\delta \sbw(c)\leq 0$ implies $\sbw\leq S_b$.
\par 2) Keeping $c=c^*$, let $\tau^{*}_o = \tau^*_{dtr} + \tau^*_{ctr} + \tau_{dpr} + \tau_{cpr}$ be perturbed to $\tau_o =\tau_{dtr} + \tau_{ctr} + \tau_{dpr} + \tau_{cpr} $ such that $\tau_{dtr}$ and $\tau_{ctr}$ satisfy 
\begin{equation}
    c^* \tau_{ctr} - \thickbar{c}^*\tau_{dtr} = \thickbar{c}^* \tau_{dpr} - c^*\tau_{cpr}, 
\end{equation}
resulting in a new bandwidth cost $\sbw(\tau_o)$. Then, $\delta \sbw (\tau_o):= \sbw-\sbw^*\leq 0$ is concave and written as:
\begin{align}
&\delta \sbw(\tau_o) = (p_3-p_4)\tau_o^2+ (q_3-q_4)\tau_o + (r_3-r_4),\label{constraintchange3} 
\end{align}
where $p_3 = -\sbw^* p_4, \ p_4 = c^*\thickbar{c}^*$,  $q_3 = 2m_{cp}n^*_{cp}\thickbar{c}^* + m_{cc} n^*_{cc}$ $c^* + 
\sbw^*(c^*\tau_{cpr}+\thickbar{c}^*\tau_{dpr})$, $q_4 = - c^* \tau_{cpr}-\thickbar{c}^*  \tau_{dpr}$, $r_3 = -\left(\sbw^* \tau_{dpr}\tau_{cpr}+ 2m_{cp}n^*_{cp}\tau_{cpr} + m_{cc}n^*_{cc}\tau_{dpr}\right), \ r_4$ $=\tau_{dpr}\tau_{cpr}$ and $\thickbar{c}^* = 1-c^*$. The constraint $\delta \sbw(\tau_o)\leq 0$ implies $\sbw\leq S_b$.} \hfill $\blacksquare$
\end{proposition}
\medskip
\par \noindent Since $\delta \sbw (\tau_o)$ and $\delta \sbw(c)$ are respectively convex and concave from Proposition \ref{propositionkcktau}, we can easily incorporate them in the co-design SDPs of Theorems \ref{theoremktau} and \ref{theoremkc} to satisfy the bandwidth constraint in \eqref{constraintchange}. Note that since $\bbK$ is co-designed with either $\tau_o$ or $c$, the true bandwidth cost $\sbw$ depends on $\bbK$ as well through $\nrow(\bbK)$, $\ncol(\bbK)$ and $\noff(\bbK)$. If $\nrow(\bbK) + \ncol(\bbK)\leq \nrow(\bbK^*) + \ncol(\bbK^*)$  and $\bs{\mf{n}}^T\noff(\bbK)\leq \bs{\mf{n}}^T\noff(\bbK^*)$, one can easily verify that $\delta \sbw(c)\leq 0$ and $\delta \sbw(\tau_o)\leq 0$ in \eqref{constraintchange2}-\eqref{constraintchange3} hold, and the true bandwidth costs always satisfy \eqref{constraintchange}. We ensure this fact by imposing a two-loop structure in our design algorithm, as seen in the next section. We next bring together the co-design SDPs \eqref{theoremktaueq}, \eqref{theoremkceq} and bandwidth constraints \eqref{constraintchange2}, \eqref{constraintchange3} in the form of our main algorithm.

\section{Problem Setup in Two-Loop ADMM Form}
\label{sec:problemsetupinadmm}
The $\mc{H}_2$-norm $J$, in general, increases with increasing sparsity of $\bbK$ \citep{negi2020sparsity}, while the bandwidth cost $\sbw$ reduces. Due to these inherent trade-offs between the objectives and the constraints, \textbf{P1} is a prime candidate to be reformulated as a two-loop ADMM optimization. The outer-loop co-designs $(\bbK,\tau_o)$ and $(\bbK,c)$ using \eqref{theoremktaueq}-\eqref{theoremkceq} under the  bandwidth constraints \eqref{constraintchange2}-\eqref{constraintchange3}. The inner-loop, on the other hand, sparsifies $\bbK$ while minimizing $J$. We describe the inner and outer loops in Sec. \ref{subsection:innerloop} and \ref{subsection:outerloop} respectively, followed by the main algorithm in Sec. \ref{subsection:mainalgo}.
\subsection{Inner ADMM Loop}
\label{subsection:innerloop}
\noindent Throughout the inner ADMM loop, we hold both $\tau_o$ and $c$ as constants. The mathematical program of the inner loop denoted as \textbf{O2.0} is written as follows:
\begin{subequations}
\label{admmform}
\begin{align}
\text{\textbf{O2.0}} \ : \ \ & \minimize_{\bbK,\bb{F}} \ \ J(\bbK) + \gamma g(\bb{F}), \\
&\text{subject to} \ \ \ \bbK=\bb{F},
\end{align}
\end{subequations}
where $\gamma$ is a regularization parameter and $g(\bb{F})=\|\bb{W}\circ \bb{F}\|_{l_1}$ is the weighted $l_1$ norm function which is used to induce sparsity in $\bb{F}$. The weight matrix $\bb{W}$ for $g(\bb{F})$ is updated iteratively through a series of reweighting steps from the solution of the previous iteration \citep{weightedl1}:
\begin{equation}
    \bb{W}_{ij} = \frac{1}{|\bb{F}_{ij}| +\epsilon }, \ \ 0<\epsilon \ll 1. \label{wl1}
\end{equation}
The augmented Lagrangian for \textbf{O2.0} is 
\begin{equation}
\mc{L}_p = J(\bbK) + \gamma g(\bb{F}) + \text{Tr}(\bs{\Theta}^T (\bbK-\bb{F})) + \frac{\rho}{2}\|\bbK-\bb{F}\|^2_{\te{F}},\label{auglag}
\end{equation}
where $\rho$ is a positive scalar, $\bs{\Theta}$ is the dual variable and $\|\cdot\|_\text{F}$ is the Frobenius norm. ADMM involves solving each objective separately while simultaneously projecting onto the solution set of the other. As shown in \cite{mihailo,boyd}, \eqref{auglag} is used to derive a sequence of iterative steps $\bbK$-min, $\bb{F}$-min and $\bs{\bs{\Theta}}$-min by completing the squares with respect to each variable.
\begin{subequations}
\label{admmeq}
\begin{align}
& \bbK_{k+1} = \underset{\bbK}{\te{argmin}} \ \Phi_1(\bbK)=\underset{\bbK}{\te{argmin}} \ J(\bbK) + \frac{\rho}{2} \|\bbK-\bb{U}_k\|^2_{\te{F}}, \label{Kmin} \\
& \bb{F}_{k+1}= \underset{\bb{F}}{\te{argmin}} \ \Phi_2(\bb{F}) =\underset{\bb{F}}{\te{argmin}} \ \gamma g(\bb{F}) + \frac{\rho}{2} \|\bb{F}-\bb{V}_k\|^2_{\te{F}}, \label{Fmin}\\
& \bs{\Theta}_{k+1} = \bs{\Theta}_k + \rho(\bbK_{k+1} - \bb{F}_{k+1}),\label{Theta}
\end{align}
\end{subequations}
where $\bb{U}_k = \bb{F}_k - \frac{1}{\rho} \bs{\Theta}_k $ and $\bb{V}_k = \bbK_{k+1}+\frac{1}{\rho}\bs{\Theta}_k$. We next present a method to solve $\bbK$-min and provide an analytical expression for $\bb{F}$-min.
\subsubsection{\texorpdfstring{$\bbK$}{K}-min Step}
Setting $\nabla \Phi_1(\bbK)=0$ and using Theorem \ref{Theoremgrad}, we get the following condition for optimality, 
\begin{equation}
\left[ (\bb{G}\bb{L}\bb{N}_d)\circ \bs{\mc{I}}_d+ (\bb{G}\bb{L}\bb{N}_o)\circ \bs{\mc{I}}_o \right] + \frac{\rho}{2}(\bbK - \bb{U}) =0, \label{kmineqgrad}
\end{equation}
where $\bb{U} = \bb{U}_k$ for the $(k+1)$-th iteration of the ADMM loop and $\bb{N}_d(c)$ is denoted as $\bb{N}_d$ as $c$ is constant for \textbf{O2.0}. $\bb{P}$ and $\bb{L}$ are the solutions of \eqref{CC} and \eqref{LL}, respectively. $\bbK$-min begins with a stabilizing $\bbK$, solves \eqref{CC}-\eqref{LL} for $\bb{P}$ and $\bb{L}$, and then solves \eqref{kmineqgrad} to obtain a new gain $\thickbar{\bbK}$ as follows: 
\begin{align}
&\thickbar{\bbK} = \texttt{Reshape}\left( (\hat{\bb{V}}_d \circ \bs{\mc{T}}_d  +  \hat{\bb{V}}_o \circ \bs{\mc{T}}_o  + \rho \bb{I}_{n^2})^{-1} \bs{\mu}   ,\texttt{[m,n]}\right),\label{tddactual}\\
&\hspace{0.25cm}\bs{\mc{T}}_d = (\bs{\mc{T}}_{dd} \circ \hat{\bb{V}}^T_d+ \bs{\mc{T}}_{od}\circ \hat{\bb{V}}^T_o), \  \bs{\mc{T}}_o = (\bs{\mc{T}}_{oo}\circ \hat{\bb{V}}^T_o + \bs{\mc{T}}_{do} \circ \hat{\bb{V}}^T_d ),  \nn\\
&\hspace{0.25cm}\bs{\mc{T}}_{dd} = 2 (\bb{N}^T_d\bb{L}\bb{N}_d \otimes \bb{R} ), \  \bs{\mc{T}}_{od} = 2(\bb{N}^T_o\bb{L}\bb{N}_d \otimes \bb{R} ), \nn \\
&\hspace{0.25cm}\bs{\mc{T}}_{oo} = 2 (\bb{N}^T_o \bb{L}\bb{N}_o \otimes \bb{R} ), \ \bs{\mc{T}}_{do} = 2(\bb{N}^T_d\bb{L}\bb{N}_o \otimes \bb{R} ),\nn \\
& \hspace{0.25cm}\bs{\mu} = \texttt{vec} \left( (2 \bs{\mc{B}}^T\bb{P}\bb{L} \bb{N}_d) \circ \bs{\mc{I}}_d +  (2 \bs{\mc{B}}^T\bb{P}\bb{L} \bb{N}_o) \circ \bs{\mc{I}}_o + \rho \bb{U} \right),\nn \\
& \hspace{0.25cm}\hat{\bb{V}}_d = \mathbf{1} \otimes \bb{v}_d, \ \bb{v}_d = \texttt{vec}(\bs{\mc{I}}_d),\ \hat{\bb{V}}_o = \mathbf{1} \otimes \bb{v}_o, \ \bb{v}_o = \texttt{vec}(\bs{\mc{I}}_o), \nn
\end{align}
For details of \eqref{tddactual}, see Appendix (Sec. \ref{subsec:derivationofkbar}). It can be shown that $\tilde{\bbK} = \bbK-\thickbar{\bbK}$ is the descent direction for $\bs{\Phi}_1$ \cite[See Lemma 4.1]{computational}. The Armijo-Goldstein line search method can then be used to determine a step size $s$ to ensure $(\bbK+s\tilde{\bbK})$ stabilizes \eqref{extendeddelayed}. The iterative process continues till we obtain $\nabla \Phi_1(\bbK)\approx 0$.

\subsubsection{\texorpdfstring{$\bb{F}$}{F}-min Step}
\label{subsub1}
The solution of the $\bb{F}$-min step is well-known in the literature \cite[Sec. 4.4.3]{boyd} as:
\begin{equation}
\bb{F}_{ij} = \begin{cases}
(1- \frac{a_{ij}}{| \bb{V}_{ij} |} ) \bb{V}_{ij}, \ \te{if} \ | \bb{V}_{ij} | > a_{ij}, \\
0, \ \te{otherwise},
\end{cases}
\end{equation}
where $a_{ij}=\frac{\gamma}{\rho} \bb{W}_{ij}$. Note that large values of $\gamma$ will induce more sparsity, and therefore may lead to a sudden increase in $J$. Therefore, $\gamma$ must be increased in small steps. { The regularization path, for example, can be logarithmically spaced from $0.01\gamma_{max}$ to $0.95\gamma_{max}$, where $\gamma_{max}$ is ideally the critical value of $\gamma$ above which the solution of \textbf{P1}\textsubscript{in} is $\bbK=\bb{F}=\bb{0}$ \citep{boyd}. In our simulations, $\gamma_{max}=1$.}
\subsection{Outer Loop}
\label{subsection:outerloop}
\noindent The outer-loop of our algorithm designs $\tau_o$ and $c$ with bandwidth constraint \eqref{constraintchange} and updates the weight matrix $\bb{W}$ for minimizing the weighted $l_1$ norm in \eqref{wl1}. Co-design of $\bbK$ in this loop is necessary to ensure stability as $\tau_o$ and $c$ change. Let $\bbK^*=\bb{F}^*$ and $\bs{\Theta}^*$ be the output of the last converged inner loop with $\bb{U}^* =\bbK^* - ({1}/{\rho})\bs{\Theta}^* $. Programs \textbf{O2.1} and \textbf{O2.2} directly design $(\bbK,\tau_o)$ and $(\bbK,c)$, respectively, in sequence as follows: 
\begin{subequations} 
\label{p101p102}
\begin{align}
\te{\textbf{O2.1}}\ : \  &\underset{\bbK,\tau_o,\bb{P}}{\te{minimize}} \ \hat{J}(\bbK,\tau_o) + \frac{\rho}{2} \|\bbK-\bb{U}^*\|^2_{\te{F}},\\
 &\te{s.t. } \ \ \ \delta \sbw(\tau_o) \leq 0, \ \text{SDP in }\te{Eq.} \  \eqref{theoremktaueq},  \\ 
\te{\textbf{O2.2}} \ :  \ &\underset{\bbK,c,\bb{L}}{\te{minimize}} \ \hat{J}(\bbK,c) + \frac{\rho}{2} \|\bbK-\bb{U}^*\|^2_{\te{F}},\\
 &\te{s.t.} \ \ \ \delta \sbw(c) \leq 0, \ \text{SDP in }\te{Eq.} \  \eqref{theoremkceq}, 
\end{align}
\end{subequations}
where $\hat{J}(\bbK,\tau_o)=\te{Tr} (\bs{\mc{B}}^T\bb{P}\bs{\mc{B}})$, $\hat{J}(\bbK,c) $ $=\te{Tr}(\bb{L} \bs{\mc{C}}^{*T}\bs{\mc{C}}^* )$, $\bs{\mc{C}}^*$ $=(\bbK^*\circ \bs{\mc{I}}_d) \bb{N}^T_d(c^*) + (\bbK^*\circ \bs{\mc{I}}_o)\bb{N}^T_o$. We next present our main algorithm to show the iterative solutions of \textbf{O2} beginning from a known stabilizing tuple $(\bbK^*,\tau^*_o,c^*)$.
\begin{algorithm}[hbtp]
\caption{Main Algorithm}
\label{MainAlgorithm}  
 \begin{algorithmic}[1]
\State \textbf{Input:} Initial feasible point $(\bbK^*_{o},\tau^*_o$, $c^*)\in\mathbb{K}$
\For{$\gamma_i = 0.01\gamma_{max}$ to $0.95\gamma_{max}$}
\State \textbf{Input:} $\bbK^*$, $\tau^*_o$ and $c^*$ stabilizing for \eqref{delayed}
\For{$1$ to Maximum Reweighted Steps}
   \State Solve \textbf{O2.1} using $\bbK^*$, $\tau^*_o$, $c^*$ to get $\hat{\bbK}$, $\tau_o$ 
  \State Solve \textbf{O2.2} using $\hat{\bbK}$, $\tau_o$, $c^*$ to get updated $\bbK$, $c$
 \State \textbf{Input:} Inner loop initial: $\bbK$, $c$, $\tau_o$ 
  \label{terminating}\While{ADMM Stopping Criteria \textbf{not} met} 
\State $\text{$\bbK${-min : }}$ Solve \eqref{Kmin} for $\bbK_{k+1}$
\State $ \text{$\bb{F}${-min : }}$ Solve \eqref{Fmin} for $\bb{F}_{k+1}$
\State Update $\bs{\Theta}$ using \eqref{Theta}
\EndWhile
\State \textbf{Result:} $\bbK^*= \bbK$, $\tau^*_o = \tau_o$, $c^* = c$
\EndFor
\State Update $\bb{W}$ using $\bbK^*$ from \eqref{wl1} and $S^*$ using $\nrow(\bbK^*)$, $\ncol(\bbK^*)$, $\bb{n}_{\bb{off}}(\bbK^*)$, $\tau^*_o$ and $c^*$ from \eqref{constraintchange}
\EndFor
\State \textbf{Result:} $\bbK$, $\tau_o$ and $c$ are obtained for $\gamma_i$
\end{algorithmic}
 \end{algorithm}
\subsection{Main Algorithm}
\label{subsection:mainalgo}
\noindent Our main algorithm is listed in Algorithm \ref{MainAlgorithm}; the following points explain its key steps.
\smallskip
\par\noindent$\bullet$ Using \textbf{O2.1}, we first co-design a stabilizing pair $(\hat{\bbK},\tau_o)$ from an initial tuple $(\bbK^*, \tau_o^*, c^*)\in\mathbb{K}$. The two are designed together as the initial $\bbK^*$ may not be stabilizing for $\tau_o$ satisfying the bandwidth constraint \eqref{constraintchange3}.
\par \noindent $\bullet$ We then use the solution of \textbf{O2.1}, i.e., $(\hat{\bbK},\tau_o,c^*)\in\mathbb{K}$ as the initial point for \textbf{O2.2} to find an updated pair $(\bbK,c)$. From Proposition \ref{propositionkcktau}, $\delta \sbw(c)$ in \eqref{constraintchange3} is convex in $c$. Let $c_{min}$ be the minimizer of $\delta \sbw(c)$. If $\hat{\bbK}$ is stabilizing for $c_{min}$, then instead of co-designing $(\bbK,c)$, we can directly set $c=c_{min}$ and $\bbK=\hat{\bbK}$, and then use a procedure similar to $\bbK$-min to minimize $J(\bbK)$, starting from $\hat{\bbK}$. The inner-loop begins with $(\bbK,\tau_o,c)\in\mathbb{K}$ and updates $\bbK$ in the direction of decreasing $J$ and increasing sparsity while { $\tau_o$ and $c$} remain constant.
\par \noindent $\bullet$ Following \citet[Sec. III-D]{mihailo} and \citet[Sec. 3.4.1]{boyd}, $\rho$ in \eqref{admmeq} is chosen to be sufficiently large to ensure the convergence of the inner ADMM loop. Since $J$ is nonconvex, convergence of this loop, in general, is not guaranteed, as is commonly seen in the sparsity promoting literature \cite{mihailo}. However, large values of $\rho$ have been shown to facilitate convergence. We use $\rho=100$ for our simulations. The stopping criterion for the inner loop in Line \ref{terminating} of Algorithm \ref{MainAlgorithm} follows \citet[Sec. 3.3.1]{boyd}. 
\begin{figure}
\centering
\includegraphics[scale=0.6]{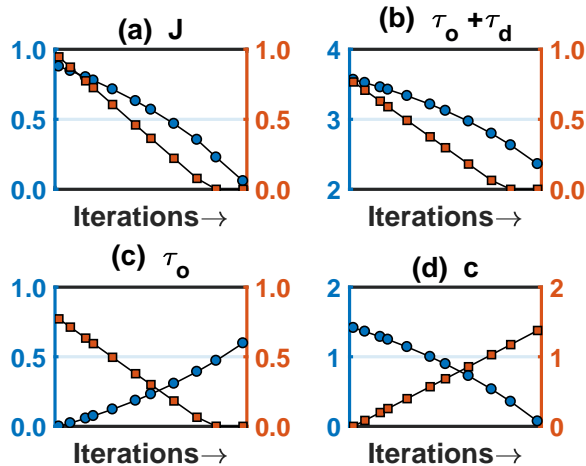}
\caption{(a), (b), (c), (d) show normalized $J$, $\tau_o+\tau_d$, $\tau_o$ and $c$ vs iterations - Right axis for {\color{orange} Model I\tsb{a}} (${\color{orange} \sqbullet}$), Left axis for {\color{cyan}Model I\tsb{b}} (${ \color{cyan}\bullet}$). }
\label{fig:case1}
\end{figure}

\section{Examples}
\label{sec:simulations1}

\subsection{Delay-Design With No Bandwidth Constraints}
\label{subsec:firstsim}

\noindent We first present simulations where only the outer loop is iterated without considering any bandwidth constraint in Algorithm \ref{MainAlgorithm}. This example shows that the relative magnitudes of $\tau_d$ and $\tau_o$ for obtaining minimum $\mc{H}_2$-norm can be significantly different for different systems. Absence of the bandwidth cost, as indicated before, will lead to the trivial solution $\tau_o=0$, $\tau_d=0$. To avoid this, we impose a simple artificial constraint $|(\tau_d-\tau^{*}_d) + (\tau_o - \tau^{*}_o)| \leq \epsilon$ where $0<\epsilon\ll 1$ is a small tolerance, and $(\bbK^*,\tau^*_o,\tau^*_d) \in  \mathbb{K}$ is the initial point for every iteration. This initial tuple is replaced by the newly designed $(\bbK,\tau_o,\tau_d)\in\mathbb{K}$ at the end of every iteration. We simulate two randomly generated models I\tsb{a} and I\tsb{b} with $\bbA\in\mathbb{R}^{5\times 5}$, $\bb{B}=\bb{B}_w=\bb{I}_n$, $\bbK^*=\bbK_{LQR}$\footnote{$\bbK_{\te{LQR}}$ is the solution of the Linear Quadratic Regulator (LQR) problem for the system in \eqref{extendeddelayed} for the given $\bb{Q}$ and $\bb{R}$.}, $\bb{Q}=\bb{R}=\bb{I}_n$ for two different initial conditions as part of Case A. The logarithm of ratios of $J$, $\tau_d+\tau_o$, $\tau_o$ and $c$ with respect to their respective minima are plotted in Fig. \ref{fig:case1}.
\smallskip
\par \noindent \textbf{Case A}: Right and left axis of all the sub-figures in Fig. \ref{fig:case1} show system I\tsb{a} and I\tsb{b} with $(c^*,\tau^*_o)$ chosen as $(0.489,0.141)$ and $(0.833,0.108)$, respectively. For both the systems, $J$ in Fig. \ref{fig:case1} (a) is seen to be decreasing as $\tau_o+\tau_d$ decreases. This is expected as $\mc{H}_2$-performance improves with a decrease in the overall delay. Fig. \ref{fig:case1} (a), (c) and (d) show that for achieving a lower $J$, the model I\tsb{a} requires a lower $\tau_o$ and a higher $c$, while I\tsb{b} requires a higher $\tau_o$ and a lower $c$. We can infer that obtaining a better $\mc{H}_2$-performance can demand completely different relative magnitudes of $\tau_d$ and $\tau_o$ depending on the system model and the initial conditions. Thus, this example validates the motivation of our problem in determining the trade-off between $\tau_d$ and $\tau_o$.
\begin{figure}[t]
\centering
\begin{minipage}{0.4\linewidth}
\centering
\includegraphics[scale=0.44]{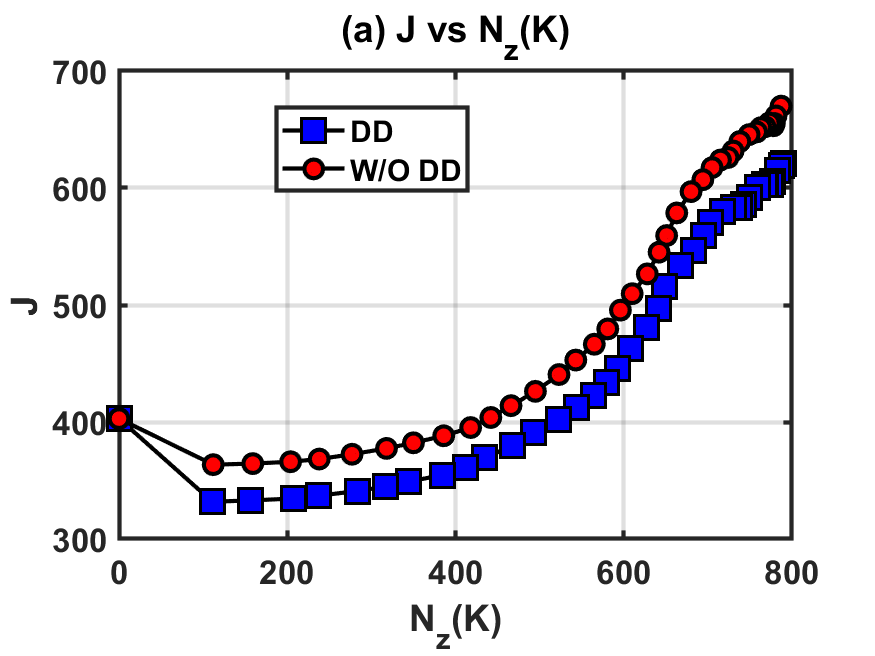} 
\end{minipage}\begin{minipage}{0.4\linewidth}
\centering
 \includegraphics[scale=0.44]{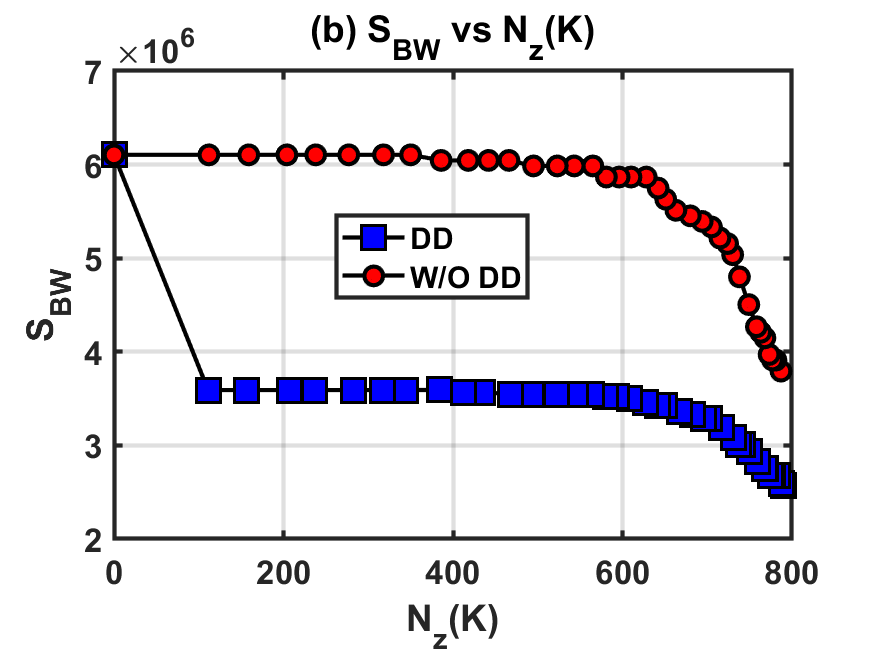}
\end{minipage}\\
\begin{minipage}{0.45\linewidth}
\centering
\includegraphics[scale=0.46]{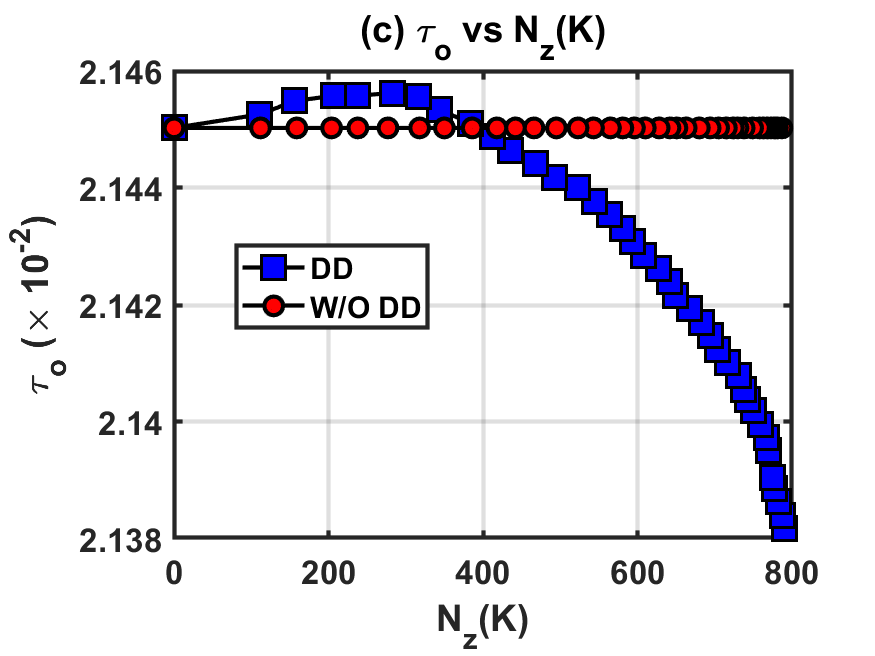} 
\end{minipage}\begin{minipage}{0.4\linewidth}
\centering
 \includegraphics[scale=0.44]{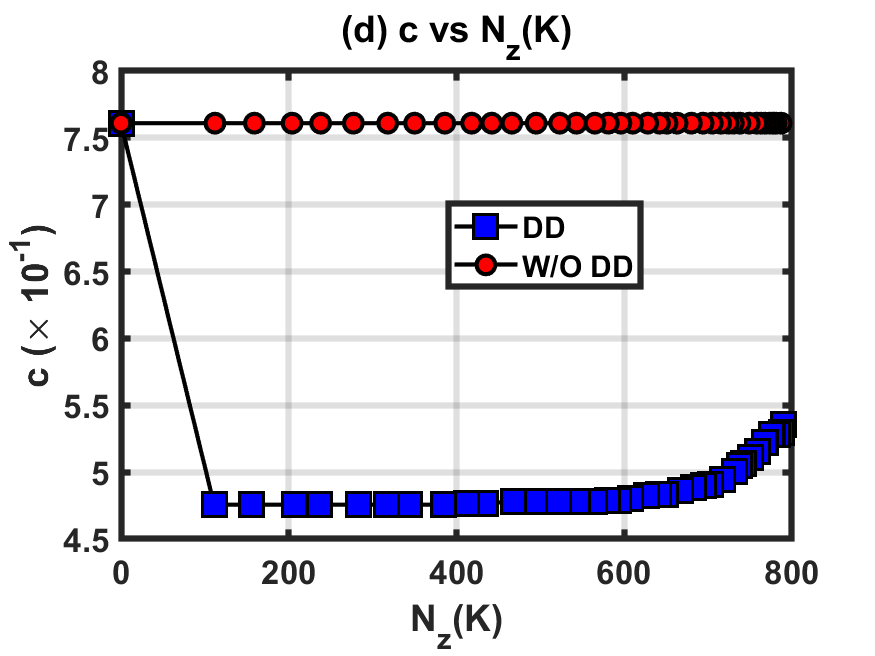}
\end{minipage}
 \caption{Case B (a), (b), (c), (d) show $J$, $\sbw$, $\tau_o$ and $c$ vs $N_z(\bbK)$ where $N_z$ is the number of zero elements of $\bbK$. `DD' and `W/O DD' indicate Algorithm \ref{MainAlgorithm} and constant-delay algorithm respectively.}
\label{fig:Case3}
\end{figure}

\subsection{Delay-Design with Bandwidth Constraints}
\label{subsec:secondsim}
\noindent We next validate Algorithm \ref{MainAlgorithm}. To illustrate its benefits, we compare it to an algorithm that consists of only the inner ADMM loop, referred to as the \textit{constant-delay} algorithm. Both algorithms start from $(\bbK^*,\tau^*_o,\tau^*_d)\in\mathbb{K}$. The delays $\tau^*_o$ and $\tau^*_d$ are kept constant throughout the constant-delay algorithm. In Case B, we present the simulations for a randomly generated LTI model ${\bbA}\in\mathbb{R}^{30\times 30}$ and $\bb{B}=\bb{B}_w=\bb{Q}=\bb{R}=\bb{I}_n$. We denote the number of zero elements of $\bbK$ by $N_z(\bbK)$.
\smallskip

\par \noindent \textbf{Case B}: We consider a randomly generated $\bbA\in\mathbb{R}^{30\times 30}$ with $900$ links in the cyber-layer, $(c^*,\tau^*_o) =(0.76,0.021)$, $m_{cp}=84$ and $m_{cc}$ $=81$. The initial conditions result in $c_{min}=0.475>c^*$ from \eqref{constraintchange2}. However, $(\bbK^*,c_{min})$ is an unstable tuple, and therefore, we rely on \textbf{O2.2} to co-design $(\bbK,c)$. Fig. \ref{fig:Case3} (a), (c) and (d) show that as sparsity increases, Algorithm \ref{MainAlgorithm} initially increases $\tau_o$ and maintains $c$ to maintain optimality of $J$. As shown in Proposition \ref{propositionkcktau}, an increase in $\tau_o$ decreases the bandwidth cost $\sbw$. Moreover, since $c$ moves towards $c_{min}$, $\sbw$ decreases steeply for Algorithm \ref{MainAlgorithm}. Further increase in sparsity of $\bbK$ causes the algorithm to decrease $\tau_o$ and increase $c$, which slows down the rate of decrease of $\sbw$ with respect to $N_z(\bbK)$. Fig. \ref{fig:Case3} (a), (b) show that as a trade-off for a much lower $\sbw$ from Algorithm \ref{MainAlgorithm}, we obtain $J$ that is comparable to the constant-delay algorithm for all the sparsity levels.

\section{Redesigning the CPS topology}
\label{sec:topologydesign}
Recall that the total system cost $S$ is composed of the \textit{bandwidth cost} $\sbw$, and the \textit{CN cost} $\scn$. In Sec. \ref{sec:problemsetupinadmm}, we designed an $\mc{H}_2$ optimal combination of sparse $\bbK$ and the delays $\tau_d$ and $\tau_c$ (i.e., bandwidths $b_{cp}$ and $b_{cc}$) to decrease $\sbw$ while keeping $\scn$ constant (i.e., keeping the topology $\bbT$ constant). In this section, as an additional step to Algorithm \ref{MainAlgorithm}, we aim to further reduce the closed-loop $\mc{H}_2$ norm of the system by redesigning the CPS topology parameters $\bs{\mathfrak{X}}$, $\bs{\mf{n}}$, $\bs{\mathfrak{U}}$ and $\bs{\mf{m}}$ while keeping the gain matrix fixed at the sparse solution $\bbK$ of Algorithm \ref{MainAlgorithm}. This redesign changes $\scn$, $\tau_d$, and $\tau_c$, which in turn changes the $\mc{H}_2$ norm $J$. In this process, $\sbw$ remains constant (i.e., the bandwidths remain fixed to the solution of Algorithm \ref{MainAlgorithm}). Since an optimal redesign of the topology should decrease both $\scn$ and $J$ beyond the solution of Algorithm \ref{MainAlgorithm}, we provide results on the existence of such an optimal topology and present a set of algorithms to obtain it.

\subsection{Effect of Topology Variation on the CPS Parameters}

Consider the system \eqref{delayed}, which is implemented as a CPS following  \ref{A1}-\ref{A4}. Before stating our design objective, we first discuss the effect of redesigning $\bbT$ with fixed $\bbK$, $b_{cc}$ and $b_{cp}$ on the following CPS parameters.
\medskip
\subsubsection{Effect on the number of Communication Links} 
\par \noindent Recall from Sec. \ref{subsubsec:sdn} and \ref{subsubsec:lan} that the number of inter-layer links is dependent on $\bbK$, while the total number of intra-layer links are dependent on both $\bbK$ and $\bb{T}$. Specifically, a non-zero $i,j$-th off-diagonal block of the matrix $\bmcK (\bbK,\bb{T})$, defined in Sec. \ref{subsubsec:sdn}, represents an intra-layer link transmitting the $n_i$ states of CN $i$ to CN $j$. 
\begin{figure}[hbtp]
    \centering
    \includegraphics[scale=0.6]{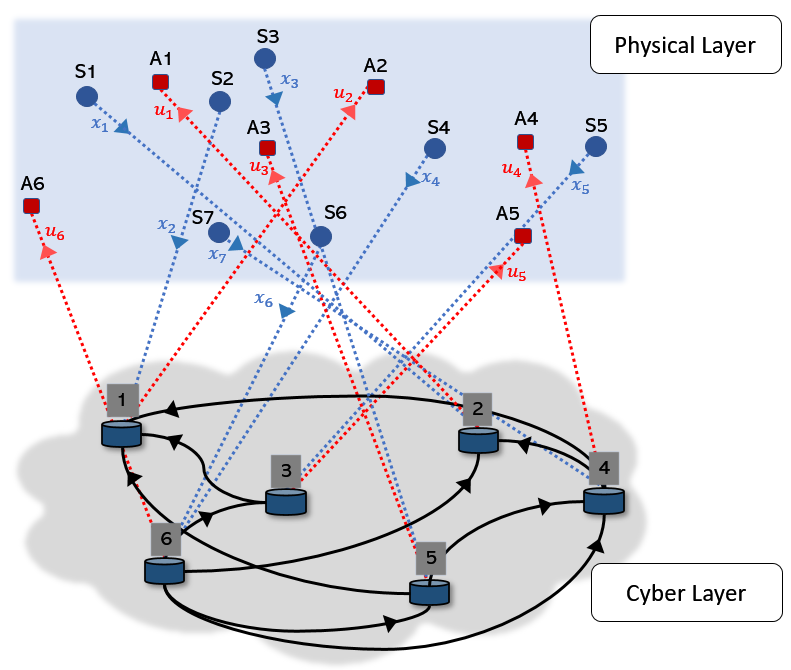}\\
    \includegraphics[scale=0.7]{sys3_newnew.png}
      \caption{A simple diagram that depicts the change in CPS topology of Fig. \ref{fig:system} that occurs due to row and column permutations of $\bbK$; See Example \ref{eg:2}. Their are $9$ intra-layer links in this topology compared to $14$ in Fig. \ref{fig:system} even though the $\nnz(\bbK)$ remains the same.}
      \label{fig:system2}
\end{figure}

\begin{example}
\label{eg:2}
{\it Let us recall Example \ref{eg:1}. For the system in Fig. \ref{fig:system}, there are a total of $\sum_i \noffs_i = 14$ intra-layer links. If we change the topology to $\bb{T}\equiv  \Big( \bs{\mf{X}} = (2,1,5,7,3,4,6)$, $\bs{\mf{U}} = (2,1,5,4,3,6)$, $\bs{\mf{n}} = [1,1,1,1,1,2]^T$, $\bs{\mf{m}} = [1,1,$ $1,1,1,1]^T \Big)$,  then, $\bmcK(\bbK,\bb{T})$ is obtained as:
  \begin{equation}
    \setlength\aboverulesep{0pt}\setlength\belowrulesep{0pt}
    \setlength\cmidrulewidth{0.5pt}
      \begin{blockarray}{cccccccc}
  {{}_{u}\mkern-1mu\setminus\mkern-1mu{}^{x}}    & {\overbrace{2}^{\bs{\mf{x}}_3}} & {\overbrace{1}^{\bs{\mf{x}}_2}} & {\overbrace{5}^{\bs{\mf{x}}_3}}& {\overbrace{7}^{\bs{\mf{x}}_4}}&
 {\overbrace{3}^{\bs{\mf{x}}_5}}&
\BAmulticolumn{2}{c}{\overbrace{4 \hspace{0.48cm} 6}^{\bs{\mf{x}}_6}} \\ 
\vspace{-0.4cm} \\
 \begin{block}{c[c|c|c|c|c|cc]}
\bs{\mf{u}}_1 \  \begin{cases} 1 \end{cases}  & e & & g  &h  &f  &  & \\
\cmidrule(lr){2-8}
\bs{\mf{u}}_2   \  \begin{cases} 2 \end{cases} &  & a &  & d  & & b& c\\
\cmidrule(lr){2-8}
\bs{\mf{u}}_3  \  \begin{cases} 3 \end{cases} &  &  & p &  &  & o & q \\
\cmidrule(lr){2-8}
\bs{\mf{u}}_4   \  \begin{cases} 4 \end{cases} &  &  &   & n  & l  & m  &  \\
\cmidrule(lr){2-8}
\bs{\mf{u}}_5  \  \begin{cases} 5 \end{cases} &  &  &  &  & i  & j & k \\
\cmidrule(lr){2-8}
\bs{\mf{u}}_6  \  \begin{cases} 6 \end{cases} &  &  &  & &  &  r & s \\
    \end{block}
\end{blockarray}. \label{Keg2}
  \end{equation} 
 Here, $\noff=[0,0,1,2,2,4]$ and therefore, the number of intra-layer links decrease from $14$ links in Example \ref{eg:1} to $\sum_i \noffs_i=9$, even though $\bbK$ remains the same. This effect of topology redesign on the assignment of the inter-layer and intra-layer links is shown in Fig. \ref{fig:system2}. As seen from the figure, the reassignment changes the destinations of the inter-layer links as well as both source and destinations of the intra-layer links. It, however, preserves the number of inter-layer links.}
\end{example}
\medskip
\subsubsection{Effect on the Delays}
\label{subsec:effectondelays}
\par { From \eqref{tauc} and \eqref{taud}, it follows that varying $\bb{T}$ for a fixed $(\bbK,$ $b_{cc},b_{cp})$ leads to a variation in $\tau_{ctr}= \kappa({\bs{\mf{n}}}^T$ $(\bbT) \noff(\bbK,\bbT))/b_{cc}$ but $\tau_{dtr} = \kappa (\nrow(\bbK) + \ncol$ $(\bbK))$ $/b_{cp})$ is kept constant (since the number of zero rows and columns will remain constant for a given $\bbK$). In Algorithm \ref{MainAlgorithm}, as $\bbT$ was fixed, both $\tau_{dpr}$ and $\tau_{cpr}$ were constant. 
In the current problem, however, $\bbT$ is the design variable, as a result of which both of these propagation delays become variable. Furthermore, we can no longer assume the per-link propagation delays for $\tau_c$ and $\tau_d$ to be equal as the choice of the destinations for the links is also variable. Therefore, we consider the worst case propagation delay for both, and accordingly replace the design variables $\tau_{cpr}$ and $\tau_{dpr}$ with $\max_i ({\tau_{cpr}}_i)$ and $\max_i ({\tau_{dpr}}_i)$ respectively, where ${\tau_{cpr}}_i$  and ${\tau_{dpr}}_i$ are the propagation delays in the corresponding $i$-th link. This will be shown shortly when we formulate our optimization problem for the redesign.}

 \subsubsection{Effect on the CN cost}
 \label{subsubsec:cncost}
\par The CN cost $\scn$ consists of two parts - the rent cost $\scnr$ and the computation overhead cost $\scnc$. 
\begin{assumption}
\label{ass:cncost}
{\it $\scnr$ is a strictly increasing function of the number of CNs $\mc{N}$. $\scnc$ is a strictly increasing function of $(n_i + m_i)$ for each CN $i$.}
\end{assumption}
Following Assumption \ref{ass:cncost}, we define the computation overhead cost as
\begin{equation}
    \scnc(\bbT)  := \sum_{i=1}^{\mc{N}} \ (n_i + m_i)^2,
\end{equation}
and the rent cost as
\begin{align}
\scnr(\bbT):= \{\scnc^{\max} \mc{X}_{(\mc{N})}| \mc{X}_{(\mc{N})} \ \text{is the $\mc{N}$-th order statistic of random samples $\mc{X}_i \sim U(0,1)$, $i\in\mathbb{N}_{\mc{N}_m}$} \},
\end{align}
where $\scnc^{\max}$ $= (m+n-2)^2 + 4$ is the maximum computation overhead cost possible for a system with $m$ inputs and $n$ states. 
\begin{figure}[t]
    \centering
    \includegraphics[scale=0.45]{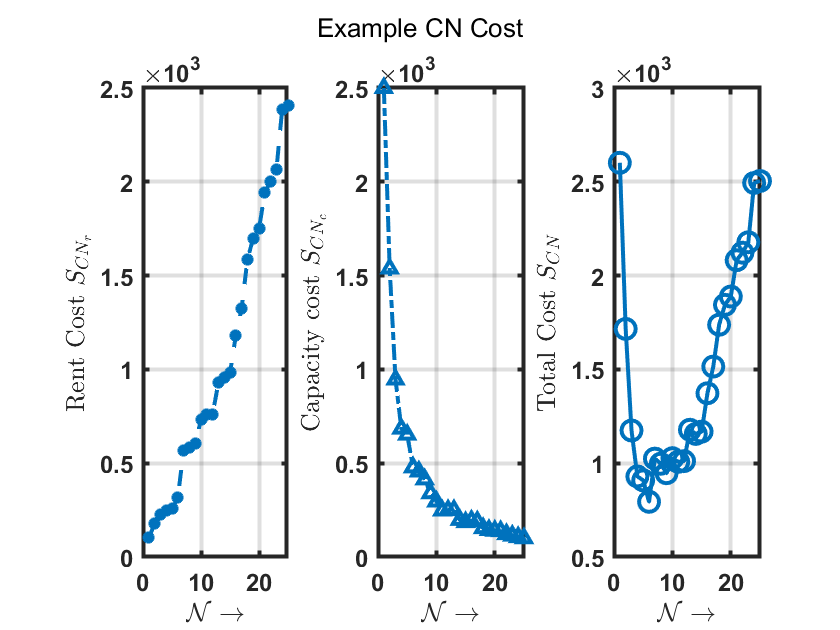}
    \caption{An example $\scn$ vs $\mc{N}$ characteristics.}
    \label{fig:costcurve}
\end{figure}
An example of total CN cost ($\scn$) vs. number of CNs ($\mc{N}$) for $m=n=25$ is shown in Fig. \ref{fig:costcurve}. We can see that the $\scn$ vs. $\mc{N}$ characteristics resemble a conventional marginal cost curve following the static pricing scheme used in cloud computing \cite{aldossary2019energy}. As the number of CNs reduces, $\scn$ first decreases due to a decrease in $\scnr$. However, as $\mc{N}$ reduces further, the increase in $\scnc$ overtakes the decrease in $\scnr$.

\subsection{Design Objective}
\par \noindent We next state our design objective as follows.
\par \textbf{P2:} {\it Let the system \eqref{delayed} be implemented as a CPS following \ref{A1}-\ref{A4} with maximum number of CNs $\mc{N}_m=\max(m,n)$ CNs (Def. \ref{def:topology}) and the initial topology denoted as $\bbT^{(in)}|_{\mc{N}_m}$. Let the outputs of Algorithm \ref{MainAlgorithm} when applied to this system be a sparse control gain matrix $\bbK^*$, bandwidths $b_{cc}^*$ and $b_{cp}^*$, and $\mc{H}_2$ norm $J(\bbK^*$, $\bbT^{(in)}|_{\mc{N}_m})$
\footnote{Since both $\tau_d$ and $\tau_c$ are functions of $\bbT$, with a slight abuse of notation, we write $J(\bbK,\bbT)$ in this section.}. Then, assigning $(\bbK,b_{cc},b_{cp})$ $\equiv$ $(\bbK^*|_{\mc{N}_m},b^*_{cc}|_{\mc{N}_m},b^*_{cp}|_{\mc{N}_m})$, find the optimal number of CNs $\mc{N}^*\in [2,\mc{N}_m]$ and the corresponding optimal topology $\bb{T}^*|_{\mc{N}^*}$ such that
\begin{enumerate}[label=\textbf{P2.\arabic*}]
    \item \label{P21} the closed-loop $\mc{H}_2$ norm is minimized and,
  \item \label{P22} the corresponding CN cost $\scn(\bbT^*|_{\mc{N}^*})$ $ \leq$ $\scn$ $(\bbT^{(in)}$ $|_{\mc{N}_m})$, i.e., it remains at most the same as the CN cost for the initial topology.
\end{enumerate} }

The rationale behind the formulation of \textbf{P2} is as follows. For a given topology $\bbT$, there is a higher loss in the closed-loop $\mc{H}_2$ norm of \eqref{delayed} when block-wise sparsity is promoted in $\bbK$ compared to when the element-wise sparsity is promoted (Algorithm \ref{MainAlgorithm}), assuming the same $\nnz(\bbK)$ in both the cases. However, block-wise sparsity in $\bbK$ is required to reduce the number of intra-layer links. Therefore, \textbf{P1} focuses on promoting element-wise sparsity in $\bbK$, and \textbf{P2} rearranges its zero entries such that a block-wise sparse structure can be obtained. When combined with the variation in $\bbT$ to change the delays, this approach enables us to reduce $J$ further than that obtained from Algorithm \ref{MainAlgorithm}. We illustrate this approach in more detail in the simulation example in Sec. \ref{subsec:blocksparse}.
\subsection{Algorithm Development}
For a fixed $\bbK^*$, each topology $\bbT|_{\mc{N}}$ of the CPS (Def. \ref{def:topology}) corresponds to an $\mc{N}\times \mc{N}$ block matrix $\bmcK|_{\mc{N}}$ through a bijective function $\bs{\mf{K}}: (\bbK^*,\bbT|_{\mc{N}} ) \to  \bmcK|_{\mc{N}} $ defined next. Let permutation matrices $\bb{U}\in \mathbb{R}^{m\times m}$ and $\bb{X}\in\mathbb{R}^{n\times n}$ be defined for a given $\bbT|_{\mc{N}}=(\bs{\mf{X}}|_{\mc{N}}$, $\bs{\mf{U}}|_{\mc{N}}$, $\bs{\mf{n}}|_{\mc{N}}$, $\bs{\mf{m}}|_{\mc{N}})$ as
 \begin{align}
    \bb{U} & :=\{\bb{I}_m : \text{ permuted as } (1,2,\ldots,m) \to \bs{\mf{U}}\}, \\
    \bb{X} & := \{\bb{I}_n : \text{ permuted as } (1,2,\ldots,n) \to \bs{\mf{X}}\}.
\end{align}
Then, the corresponding $\bmcK|_{\mc{N}}$ is obtained as:
\begin{align}
  \bmcK|_{\mc{N}} = \bs{\mf{K}}(\bbK^*,\bbT|_{\mc{N}}) =  [\bb{U} \bbK^* \bb{X}]_{[\bs{\mf{n}}|_{\mc{N}},\bs{\mf{m}}|_{\mc{N}}]}. \label{permute}
\end{align}
Here, $[\bb{X}]_{[\bb{a},\bb{b}]}$ represents the partitioning of $\bb{X}$ using rowgroups $\bb{a}$ and colgroups $\bb{b}$ (see Example \ref{eg:1}). Since $\bs{\mf{K}}$ is invertible, the set of block matrices $\bmcK|_{\mc{N}}$ that each correspond to a valid topology for $\mc{N}$ CNs can be written as:
 \begin{align}
 &\{\bmcK|_\mc{N}\} := \{\bmcK|_\mc{N} : \bbT|_\mc{N} = \bs{\mf{K}}^{-1}(\bbK^*,\bmcK|_\mc{N})  \text{ defines a topology from Def. \ref{def:topology} } \}. \label{feasibleset}
 \end{align}

 \par To solve \textbf{P2}, we first choose optimal $\bmcK^*|_{2}$, $\cdots$, $\bmcK^*|_{\mc{N}_m}$ from the sets $\{\bmcK|_{2}\}$, $\cdots$, $\{\bmcK|_{\mc{N}_m}\}$ such that their corresponding CN costs fulfill \ref{P22}; the guidelines for this choice will be defined shortly. Out of these $\mc{N}_m-1$ matrices, we choose $\bmcK^*|_{\mc{N}^*}$ as the one that corresponds to the minimum closed-loop $\mc{H}_2$ norm. The solution of \textbf{P2} will then be obtained as $\bbT^*|_{\mc{N}^*}$ $=$ $\bs{\mf{K}}^{-1}(\bbK^*,\bmcK^*|_{\mc{N}^*})$. We begin with the steps for obtaining $\{\bmcK|_{\mc{N}}\}$ in \eqref{feasibleset} for $\mc{N}=2$.

\subsubsection{Block Diagonal Permutation}
\label{subsubsec:rowcolpermute}
To obtain $\bmcK^*|_{\mc{N}=2}$, we first need to find $\bb{U}$ and $\bb{X}$ such that the majority of the zeros in $\bbK^*$ are delegated to the off-diagonal blocks $\bmcK_{1,2}$ and $\bmcK_{2,1}$.
\begin{equation}
    \bmcK|_2 = [\bb{U} \bbK^* \bb{X}]_{[\bs{\mf{n}}|_2,\bs{\mf{m}}|_2]} = \begin{bmatrix}
    \bmcK_{1,1} & \bmcK_{1,2} \\ 
    \bmcK_{2,1} & \bmcK_{2,2}
    \end{bmatrix}. \label{partion1}
\end{equation}
This is done so that the total number of non-zero off-diagonal blocks $\sum \noff(\bmcK|_{\mc{N}=2})$, i.e., the number of intra-layer links decreases. For a fixed $b^*_{cc}$, this can result in a decrease in the intra-layer transmission delay $\tau_{ctr}$ and eventually in the closed-loop $\mc{H}_2$ norm.
\par  We use the spectral partitioning method based on the Fiedler vector to carry out the block diagonal permutations, following its variant for rectangular matrices as presented in \cite{kolda}. Algorithm \ref{algo:spectralpartitioning} states the steps for carrying out these permutations. Since the block sizes, i.e., $\bs{\mf{n}}|_2$ and $\bs{\mf{m}}|_2$ are unknown, we assign $m_1= \texttt{floor}[ m/2]$ and $n_1 = \texttt{floor}[n/2]$ to avoid a trivial solution. The actual block sizes are determined after obtaining $\bb{U}\bbK^*\bb{X}$, as seen next.

\subsubsection{Optimal Partitioning}
\label{subsubsec:optimalpart}
Once $\bb{U}$ and $\bb{X}$ are obtained from Algorithm \ref{algo:spectralpartitioning}, we next search for sets of $\bs{\mf{n}}|_2$ and $\bs{\mf{m}}|_2$ to obtain the feasible set $\{\bmcK|_2\}$. There are $(m-1)\times(n-1)$ ways to partition $\bb{U}\bbK^* \bb{X}$ as a $2 \times 2$ block matrix. For example, if $m=n=3$, then $\bmcK$ could be partitioned into a $2\times 2$ block matrix in the following $4$ ways:
\begin{align}
\begin{array}{cccccc}
& \left[\begin{array}{c|cc}
\bullet & \sqbullet & \bullet\\
\hline
\bullet & \varstar & \blacktriangleup\\
\blacktriangleup & \bullet & \blackdiamond
\end{array}\right]  &       
\left[\begin{array}{c|cc}
\bullet & \sqbullet & \bullet\\
\bullet & \varstar & \blacktriangleup\\
\hline
\blacktriangleup & \bullet & \blackdiamond
\end{array}\right] & 
\left[\begin{array}{cc|c}
\bullet & \sqbullet & \bullet\\
\hline
\bullet & \varstar & \blacktriangleup\\
\blacktriangleup & \bullet & \blackdiamond
\end{array}\right]  &  \left[\begin{array}{cc|c}
\bullet & \sqbullet & \bullet\\
\bullet & \varstar & \blacktriangleup\\
\hline
\blacktriangleup & \bullet & \blackdiamond
\end{array}\right] \\[1cm]
\bs{\mf{n}}|_2 & [1,2] &  [1,2]  &  [2,1]  & [2,1]\\[0.5cm]
\bs{\mf{m}}|_2  & [1,2] & [2,1] & [1,2] & [2,1] 
\end{array}
\end{align}
\refstepcounter{alphaforalgorithm}
\renewcommand{\thealgorithm}{\arabic{algorithm}\alph{alphaforalgorithm}}
\begin{algorithm}[t]
\hypertarget{algo:2}{}
\caption{Spectral Partitioning Algorithm}
\label{algo:spectralpartitioning} \begin{algorithmic}[1]
\State \textbf{Input:} Sparse matrix $\bbK^*\in\mathbb{R}^{m\times n}$ 
\State Symmetrize $\bbK^*$ as $\hat{\bbK}^* = [ \bb{0} \ \bbK^* ;  {\bbK^*}^T  \ \bb{0}] $
\State Compute Laplacian of $\hat{\bbK}^*$ : $\bb{L} = \bb{D}- \hat{\bbK}^*$, where $\bb{D}=\texttt{Diag}\{d_1,d_2,\ldots,d_{m+n}\}$ and $d_i = \sum_j \hat{\bbK}^*_{ij}$ 
\State Find $\texttt{w}$ as the eigenvector corresponding to the smallest positive eigenvalue of $\bb{L}$ or Fiedler vector
\State $\texttt{w}_{f}\in\mathbb{R}^{m}$ and $\texttt{w}_{l}\in\mathbb{R}^n$ denote first $m$ and last $n$ elements of $\texttt{w}$, respectively. 
\State Sort $\texttt{w}_{f}$ and $\texttt{w}_{l}$ in descending order and denote the corresponding sorted indices as $\{i_1,i_2,\ldots,i_m\}$ and $\{j_1,j_2,\ldots,j_n\}$. 
\State $\bb{U} :=\text{$\bb{I}_m$ with rows permuted as $\{i_1,i_2,\ldots,i_m\}$ }$
\State $\bb{X} :=\text{$\bb{I}_n$ with rows permuted as $\{j_1,j_2,\ldots,j_m\}$ }$
\State \textbf{Result:} Permutation matrices $\bb{U}$, $\bb{X}$, rowgroup $\bs{\mf{U}}=\{i_1,i_2,\ldots,i_m\}$ and colgroup $\bs{\mf{X}}=\{j_1,j_2,\ldots,j_m\}$
\end{algorithmic}
 \end{algorithm}
  \setcounter{algorithm}{1}
\refstepcounter{alphaforalgorithm}
\renewcommand{\thealgorithm}{\arabic{algorithm}\alph{alphaforalgorithm}}
\begin{algorithm}[t]
\caption{Optimal Partitioning Algorithm}
\label{algo:optimalpartitioning}
\begin{algorithmic}[1]
\State \textbf{Input:} Initial sparse gain matrix $\bbK^*$, CN Cost $\scn(\bbT^{(in)}|_{\mc{N}_m})$, $\mc{H}_2$ norm $J(\bbK^*,\bbT^{(in)}|_{\mc{N}_m})$
\State \textbf{Input:} $[\bb{U}\bb{K^*}\bb{X}]$, $\bs{\mf{U}}$ and $\bs{\mf{X}}$ from Algorithm \ref{algo:spectralpartitioning}
\For {$i=1$ till $m-1$}
  \State Set $\bs{\mf{m}} =[i,m-i]^T$ 
 \For{$j=1$ till $n-1$} 
   \State Set $\bs{\mf{n}}=[j,n-j]^T$ 
   \State Set $\bbT := (\bs{\mf{X}},\bs{\mf{U}},\bs{\mf{n}},\bs{\mf{m}}) $ for current iteration
   \State Calculate $\tau_c(\bbT)$, $\tau_d(\bbT)$ and $S_{CN}(\bbT)$
   \State Calculate $J(\bbK^*,\tau_c,\tau_d)$ for the current iter
   \Statex \hspace{1cm}-ation's $\bs{\mf{m}}$ and $\bs{\mf{n}}$
   \EndFor
   \EndFor
\State Choose the partitionings $\bs{\mf{m}}$, $\bs{\mf{n}}$ for which the corresponding $S_{CN}(\bbT)$ $\leq$ $S_{CN}(\bbT^{(in)}|_{\mc{N}_m})$ 
\State Out of the remaining partitionings, choose the one that corresponds to the lowest $\mc{H}_2$ norm and set $\bbT^*|_{\mc{N}} = \bbT$ for that partitioning
\State \textbf{Result:} Topology $\bbT^*|_{\mc{N}}$ for $\mc{N}$ number of CNs
\end{algorithmic}
 \end{algorithm}
\par \noindent The computation of $\{\bmcK|_2\}$ is detailed in Algorithm \ref{algo:optimalpartitioning}. We next discuss the variation of $\tau_c$, $\tau_d$, $\scn$ over the set $\{\bmcK|_\mc{N}\}$.
\par $\bullet$ \textbf{$\tau_d(\bmcK|_\mc{N})$:} While $\tau_{dtr}(\bbK^*)$ remains constant, $\tau_{dpr}(\bbT$ $|_\mc{N})$ $\equiv \tau_{dpr}(\bs{\mf{K}}^{-1}(\bbK^*,\bmcK|_\mc{N}))$\footnote{\label{footnote:inversefunction}Here, $\bs{\mf{K}}^{-1}(\bbK^*,\bmcK|_{\mc{N}})$ is the inverse of the function $\bs{\mf{K}}(\bbK^*$ $\bbT|_\mc{N})$ in the argument $\bbT|_{\mc{N}}$ for fixed $\bbK^*$, as given in \eqref{permute}.} varies over $\{\bmcK|_{\mc{N}}\}$.
\par $\bullet$ \textbf{$\tau_c(\bmcK|_\mc{N})$:} Both $\tau_{ctr}$ and $\tau_{cpr}$ are functions of $\bbT$ and therefore, can vary over $\{\bmcK|_\mc{N}\}$. To ensure that $\tau_{cpr}$ remains constant or decreases over $\{\bmcK|_{\mc{N}}\}$ for $\mc{N} <\mc{N}_m$, we can choose an $\mc{N}$-dimensional set of CNs (out of $\mc{N}_m$). We formally state the existence of such a set in the following lemma.
\smallskip
    \begin{lemma}
    \label{lemma:taucpr1}
  {\it Assume a CPS with parameters $(\bbK^*$, $b^*_{cc}$, $\mc{N}_m)$ and topology $\bbT^{(in)}|_{\mc{N}_m}$ that results in the intra-layer propagation delay $\tau_{cpr}(\bbT^{(in)}|_{\mc{N}_m})$. For any $\mc{N}\leq \mc{N}_m$, one can always choose an $\mc{N}$-dimensional set of CNs such that for all topologies $\bbT|_{\mc{N}}$, $\tau_{cpr}(\bbT|_{\mc{N}})\leq \tau_{cpr}(\bbT^{(in)}|_{\mc{N}_m})$.} \hfill $\blacksquare$
    \end{lemma}
\smallskip
\par $\bullet$ \textbf{$\scn(\bmcK|_\mc{N})$:} Rent cost $\scnr$ is constant but the computation overhead cost $\scnc(\bbT|_\mc{N})\equiv \scnc$ $(\bs{\mf{K}}^{-1}(\bmcK|_{\mc{N}}))$\footref{footnote:inversefunction} varies over $\{\bmcK|_\mc{N}\}$.
 \par Our objective is to find the optimal $\bmcK^*|_\mc{N}$ out of $\{\bmcK|_\mc{N}\}$ that minimizes the closed-loop $\mc{H}_2$ norm under the constraint that the CN cost does not increase from that in Algorithm \ref{MainAlgorithm}.
\begin{align}
\textbf{O3:} \   \bmcK^*|_\mc{N} =  \ &\underset{\{\bmcK|_\mc{N}\}}{\text{argmin}} \ J(\tau_c(\bmcK|_\mc{N}),\tau_d(\bbT(\bmcK|_\mc{N})),\bbK^* ), \\
    &\text{s.t.} \ \ \  \scn(\bmcK|_\mc{N}) \leq \scn(\bbT^{(in)}|_{\mc{N}_m}).
\end{align}
 Algorithm \hyperlink{algo:2}{2} obtains $\bbT^*|_{{2}} \equiv \bs{\mf{K}}^{-1}(\bbK^*,\bmcK^*|_2)$\footref{footnote:inversefunction} by solving \textbf{O3}. This is the optimal topology that can be implemented for a $2$ CN CPS. To obtain the same for $\mc{N}>2$ CNs, the diagonal blocks of $\bmcK|_2$ are further divided into their constituent $2\times 2$ blocks recursively by using Algorithm \hyperlink{algo:2}{2}. Once all the diagonal blocks are exhausted, i.e., cannot be further divided, the off-diagonals blocks follow the same recursive division process, until we obtain $\bbT^*|_{\mc{N}_m}\equiv \bs{\mf{K}}^{-1}(\bmcK^*|_{\mc{N}_m})$. Finally, the optimal topology from the designed $\bbT^*|_2$, $\bbT^*|_3$, $\cdots$, $\bbT^*|_{\mc{N}_m}$ is chosen such that $\min_\mc{N} J(\bbT^*|_\mc{N})$. The optimization \textbf{O3} is trivially feasible for $\mc{N}_m$ CNs. We next state the proposition that derives the conditions under which \textbf{O3} is feasible for any $\mc{N}<\mc{N}_m$.

\begin{proposition}
\label{prop:sec6}
{\it Assume a CPS with $\mc{N}_m$ CNs and fixed $(\bbK^*$, $b^*_{cc}$, $b^*_{cp})$ has a topology $\bbT^{(in)}|_{\mc{N}_m}$, intra-link transmission delay $\tau_{ctr}(\bbT^{(in)}|_{\mc{N}_m})$, and CN cost $\scn(\bbT^{(in)}$ $|_{\mc{N}_m})$. Let Algorithm \hyperlink{algo:2}{2} be used to obtain a topology $\bbT^*|_{\mc{N}}$ for $\mc{N}<\mc{N}_m$ such that the fraction of the total non-zero off-diagonal blocks in $\bmcK(\bbK^*,\bbT^*|_{\mc{N}})$ is less than or equal to that in $\bmcK(\bbK^*,\bbT^*|_{\mc{N}_m})$, i.e.,
\begin{equation}
\frac{\sum_i^\mc{N} {n_{\text{off}}}_i(\bbK^*,\bbT^*|_{\mc{N}})}{\mc{N}(\mc{N}-1)} \leq \frac{\sum_i^{\mc{N}_m} {n_{\text{off}}}_i(\bbK^*,\bbT^{(in)}|_{\mc{N}_m})}{\mc{N}_m(\mc{N}_m-1)}. \label{condition1}
\end{equation}
When $n\leq m$, there always exists an optimal number of CNs $\mc{N}^*<\mc{N}_m$ for which $\tau_{ctr}(\bbT^*|_{\mc{N}^*})< \tau_{ctr}(\bbT^{(in)}$ $|_{\mc{N}_m})$ and $\scn(\bbT^*|_{\mc{N}^*})$ $\leq$ $\scn(\bbT^{(in)}|_{\mc{N}_m})$. When $n>m$. then this result holds if $n \leq \frac{1}{2} m(m-1)$.}  \hfill $\blacksquare$
\end{proposition}
\par \noindent Lemma \ref{lemma:taucpr1} and Proposition \ref{prop:sec6} provide sufficient conditions under which intra-layer delay $\tau_c(\bbT^*|_{\mc{N}^*})$ $<\tau_c(\bbT^{(in)}$ $|_{\mc{N}_m})$ for some $\mc{N}^*$ $<\mc{N}_m$. If the increase in the corresponding $\tau_{d}(\bbT^*|_{\mc{N}^*})$ is limited such that $(\tau_c(\bbT^*|_{\mc{N}^*})+\tau_d(\bbT^*|_{\mc{N}^*}))$ $<$ $(\tau_c(\bbT^{(in)}|_{\mc{N}_m})+\tau_d(\bbT^{(in)}|_{\mc{N}_m}))$, one can obtain $J(\bbT^*|_{\mc{N}^*})< J(\bbT^{(in)}|_{\mc{N}_m})$. Otherwise, one should stick to the initial topology $\bbT^{(in)}|_{\mc{N}_m}$. Note that $\tau_{cpr}(\bbT^*|_{\mc{N}^*})$ may also be lesser than the initial value $\tau_{cpr}(\bbT^{(in)}|_{\mc{N}_m})$, as shown in Lemma \ref{lemma:taucpr1}. This can result in a further decrease in the overall delay $\tau_c(\bbT^*|_{\mc{N}^*}) + \tau_d(\bbT^*|_{\mc{N}^*})$.
\begin{figure}[h!]
    \centering
      \includegraphics[scale=0.35]{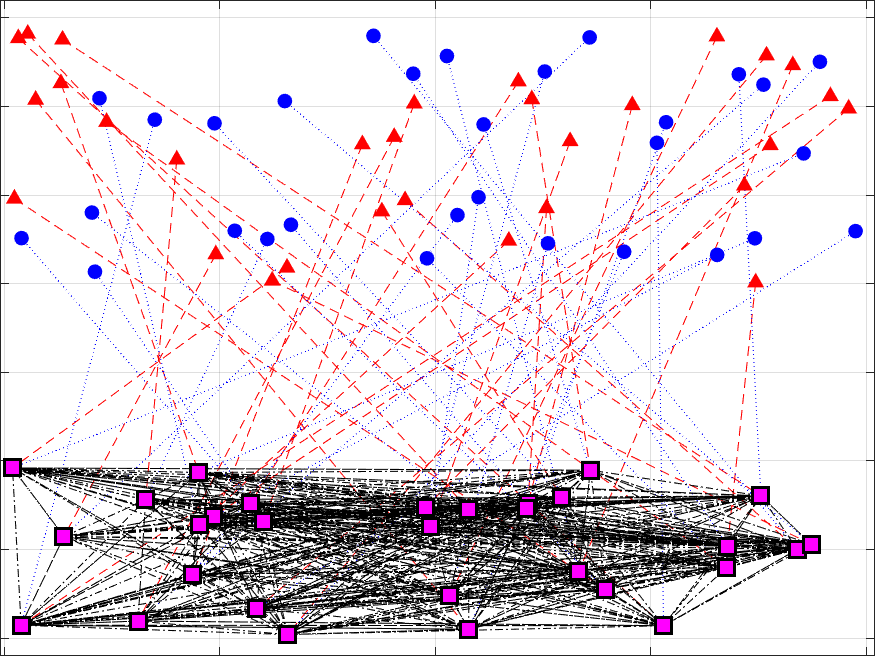}
     \includegraphics[scale=0.35]{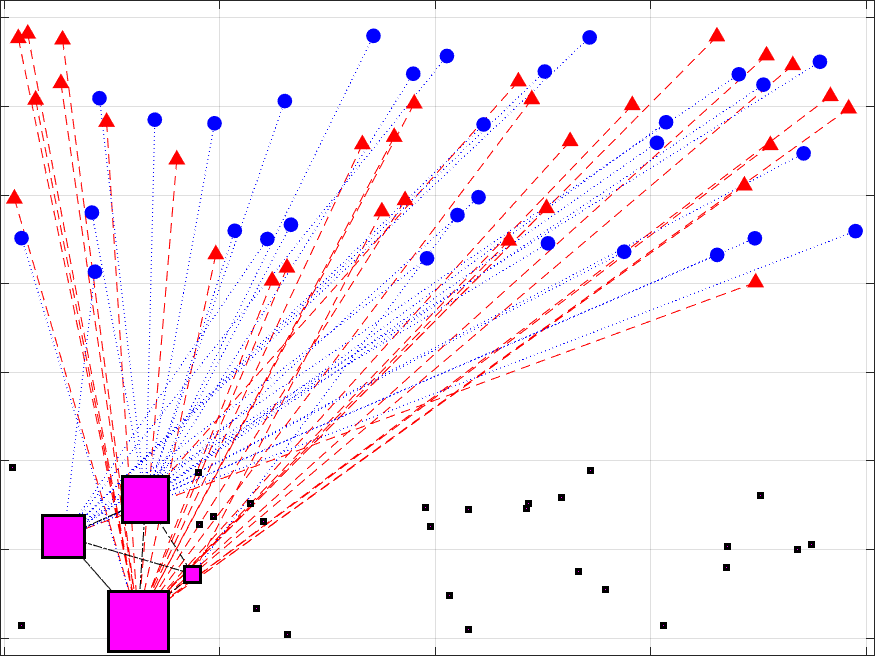} 
    \includegraphics[scale=0.35]{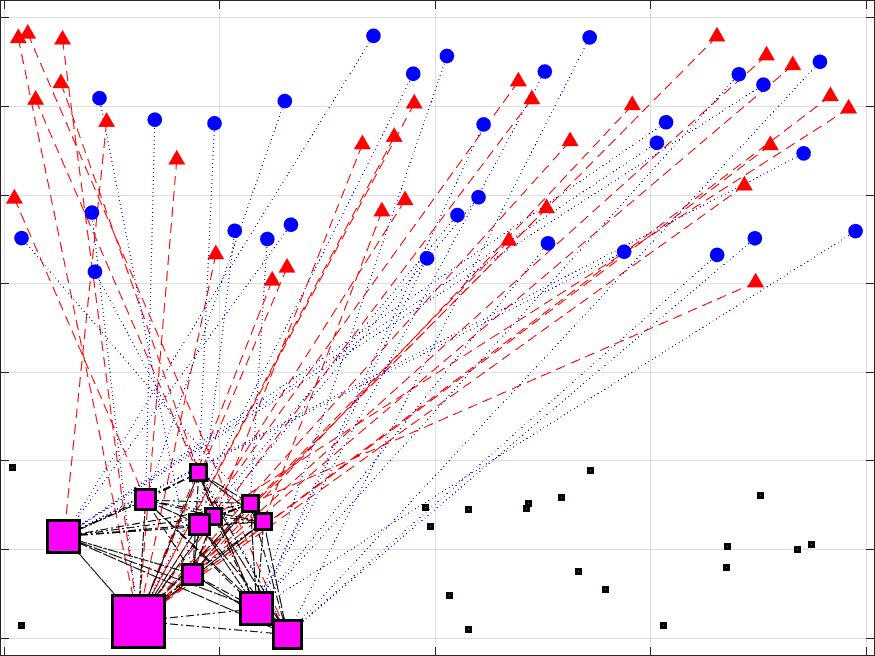}
    \\
      (a) $\bbT^{(in)}|_{\mc{N}_m=30}$
\hspace{3cm}  (b) $\bbT^*|_{\mc{N}=4}$ \hspace{4cm} (c) $\bbT^*|_{\mc{N}=11}$\\
\vspace{0.5cm}
    \includegraphics[scale=0.35]{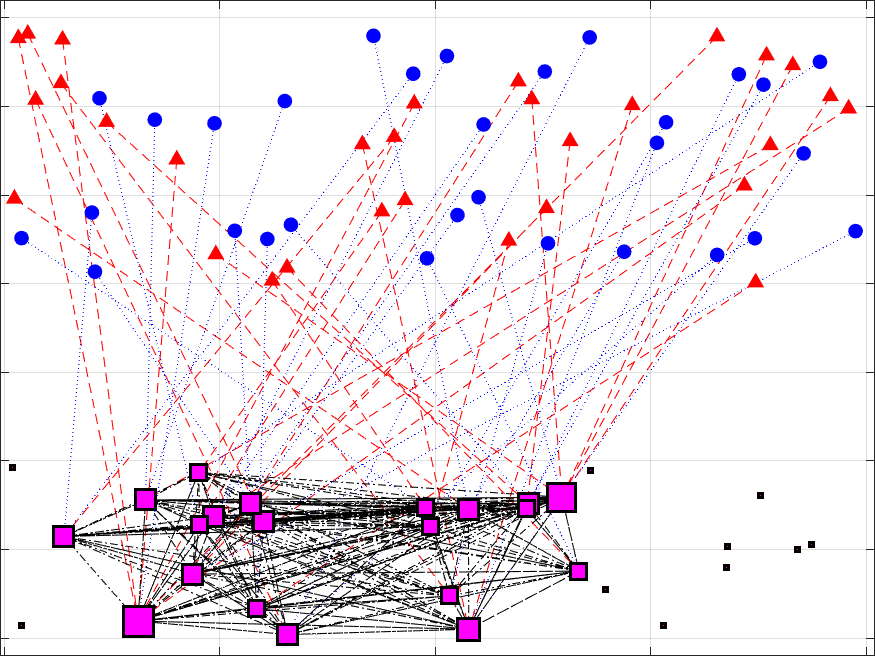} 
     \includegraphics[scale=0.35]{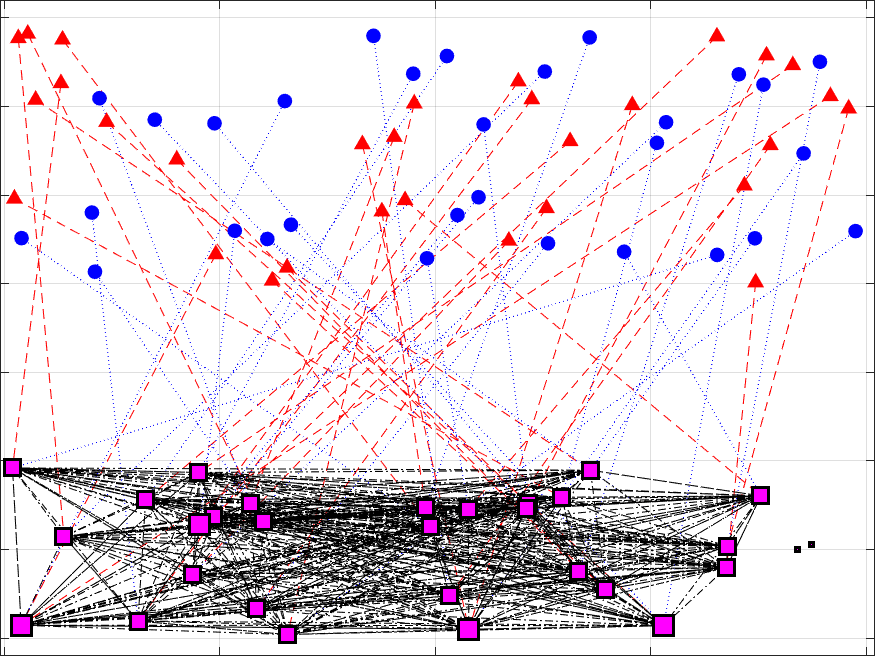}
         \includegraphics[scale=0.35]{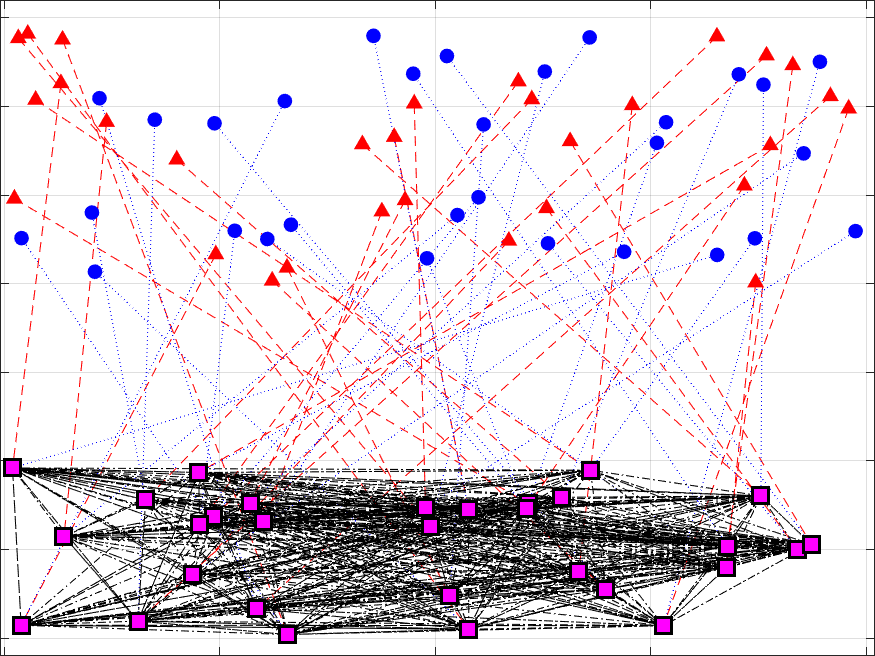}\\
     (d) $\bbT^*|_{\mc{N}_m=20}$
\hspace{3.2cm}  (e) $\bbT^*|_{\mc{N}=28}$ \hspace{4cm} (f) $\bbT^*|_{\mc{N}=30}$\\
  \begin{center}  {\red  \large $\blacktriangleup$} Actuator \hspace{1.5cm}  { \large $\bullet$} Sensor \hspace{1.5cm} {\color{magenta} \large $\sqbullet$} CN \hspace{1.5cm} { \color{red} $--$} Inter-layer link carrying $\bbu(t)$  \end{center}  
 \begin{center}
        {\color{blue} $\cdots$} Inter-layer link carrying $\bbx(t)$ \hspace{0.5cm} {\color{black} $-\cdot-$} Intra-layer link carrying $\bbx(t)$ 
        \end{center}
    \caption{Case B: (a) shows the initial topology $\bbT^{(in)}|_{\mc{N}_m}$ for $\mc{N}_m=30$ CNs used in Algorithm \ref{MainAlgorithm}. The optimal topologies for $\mc{N} = 4$, $11$, $20$, $28$ and $30$ obtained from Algorithm \texorpdfstring{\protect\hyperlink{algo:2}{2}}{2} are shown in (b), (c), (d), (e) and (f) respectively.}
        \label{fig:caseb2CLcombplot}
\end{figure}
\section{Examples}
\label{sec:simulations2}
We apply Algorithm \hyperlink{algo:2}{2} to the simulation example of Case B presented in Sec. \ref{sec:simulations1}. A sparse $\bbK^*$ obtained from Algorithm \ref{MainAlgorithm} with the corresponding topology $\bbT^{(in)}|_{\mc{N}_m=30}$ is used as the input to Algorithm \hyperlink{algo:2}{2}. The values of $\bs{\mathfrak{X}}^{(in)}|_{30}$, $\bs{\mathfrak{U}}^{(in)}|_{30}$, $\bs{\mf{n}}^{(in)}|_{30}$ and $\bs{\mf{m}}^{(in)}|_{30}$ for this initial topology are randomly chosen and given in {\red negi (2021)}. Moreover, the values of $\tau_d(\bbT^{(in)}|_{30})$, $\tau_c(\bbT^{(in)}|_{30})$, $\scn(\bbT^{(in)}|_{30})$ and $J(\bbK^*,$ $\bbT^{(in)}|_{30})$ are used to normalize the corresponding $\tau_d(\bbT^*|_{\mc{N}})$, $\tau_c(\bbT^*|_{\mc{N}})$, $\scn(\bbT^*|_{\mc{N}})$ and $J(\bbK^*,\bbT^*|_{\mc{N}})$ for optimal topologies $\bbT^*|_{\mc{N}}$ obtained as the outputs of Algorithm \hyperlink{algo:2}{2}. For e.g., a normalized $\tau_c(\bbT^*|_{\mc{N}})=0.8$ means that for some $\mc{N}<\mc{N}_m$, Algorithm \hyperlink{algo:2}{2} obtains a topology with $20\%$ less intra-layer delay compared to the initial topology.

\subsection{Applying Algorithm \texorpdfstring{\protect\hyperlink{algo:2}{2}}{2} to Case B}
\begin{figure}[h!]
    \centering
    \includegraphics[scale=0.65]{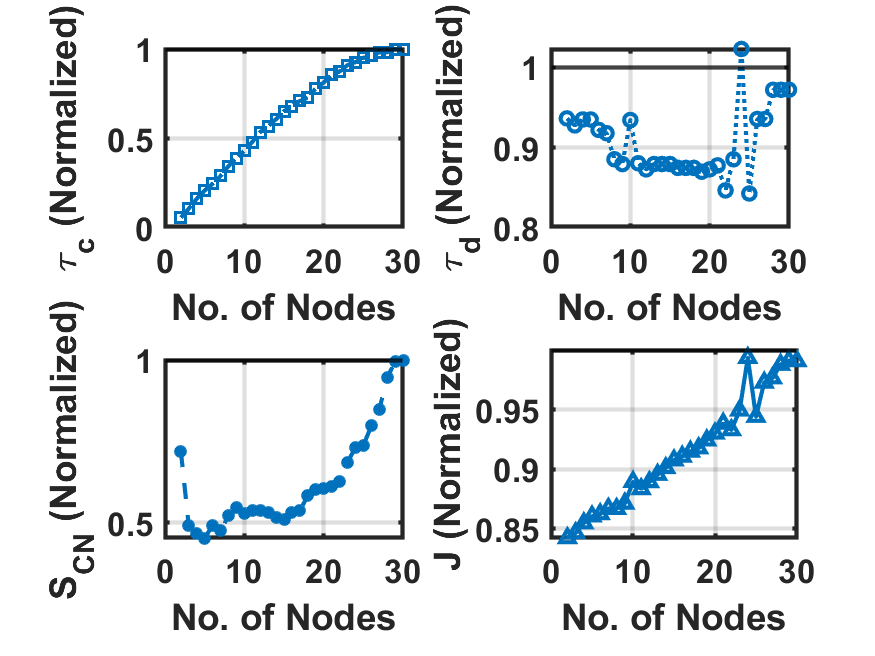}
   \caption{Case B: Comparison of normalized $J$, $S$, $\tau_d$ and $\tau_o$ between different designed topologies $\bbT^*|_2$, $\bbT^*|_3$, $\cdots$ $\bbT^*|_{30}$. The normalization is with respect to the corresponding parameter values of the initially given topology.}
    \label{fig:caseb2combToplot}
\end{figure}

{ The input to Algorithm \hyperlink{algo:2}{2}
 is a sparse $\bbK^*\in\mathbb{R}^{30 \times 30}$ with $\nnz=131$. Fig. \ref{fig:caseb2CLcombplot} (a) shows the initial topology $\bbT^{(in)}|_{30}$, while Fig. \ref{fig:caseb2CLcombplot} (b) - (f) show the optimal topologies $\bbT^*|_{\mc{N}}$ obtained as the output of Algorithm \hyperlink{algo:2}{2} for $\mc{N} = 4$, $11$, $20$, $28$ and $30$, respectively. In these figures, the size of the square representing a CN is directly proportional to the number of states and inputs associated with it. Thus, these figures show how Algorithm \hyperlink{algo:2}{2} modifies the number of states and inputs associated with each CN to satisfy \textbf{O3} as $\mc{N}$ changes. Furthermore, they show the subset of $\mc{N}$ chosen out of $\mc{N}_m$ CNs to fulfill Lemma \ref{lemma:taucpr1}. Fig. \ref{fig:caseb2combToplot} shows that Algorithm \hyperlink{algo:2}{2} guarantees a lower value of the closed-loop $\mc{H}_2$ norm for all $\mc{N}\leq \mc{N}_m$ compared to $J(\bbK^*,\bbT^{(in)}|_{\mc{N}_m})$. For the optimal topologies, any potential increase in $\tau_{dpr}$ is compensated by the decrease in $\tau_{ctr}$, as seen for $\mc{N}=24$. The CN cost $\scn$ decreases continuously from $\mc{N}=30$ until $\mc{N}=15$, after which it starts increasing again as the sharp rise in $\scnc$ overshadows the decrease in the cost of renting a lower number of CNs. The main message conveyed by Fig. \ref{fig:caseb2combToplot} is that  if one desires the lowest $J$ such that $\scn\leq \scn(\bbT^{(in)}|_{\mc{N}_m}) $, then one can choose the optimal topology $\bbT^{*}|_{\mc{N}}$ obtained for $\mc{N}=2$; however, if the goal is to obtain the lowest $\scn$ such that $J\leq J(\bbT^{(in)}|_{\mc{N}_m})$, then one must choose the optimal topology corresponding to $\mc{N}=4$.}
\subsection{Algorithm \texorpdfstring{\protect\hyperlink{algo:2}{2}}{2} vs Block-Sparse Algorithm of \texorpdfstring{\citet{negi2020sparsity}}{[]}}
\label{subsec:blocksparse}
\par
\begin{figure}[h]
    \centering
    \includegraphics[scale=0.4]{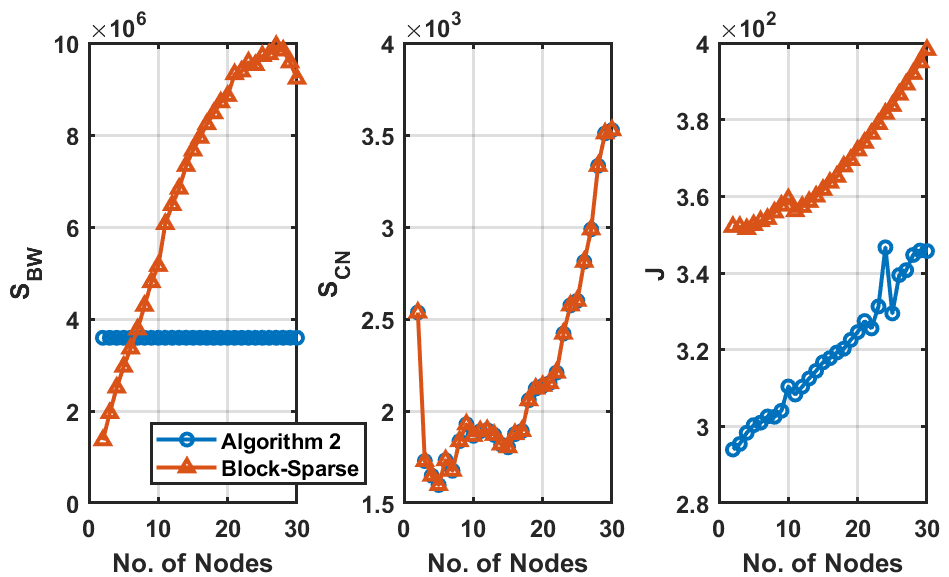} 
    \caption{Comparison of $\sbw$, $\scn$ and $J$ between Algorithm \texorpdfstring{\protect\hyperlink{algo:2}{2}}{2} output and direct block-sparsity promotion in $\bbK^*$ for Case B}
    \label{fig:caseb2BlockSparsePlot}
\end{figure}

 Note that our approach is to first promote element-wise sparsity in the initial $\bbK^*$ using Algorithm \ref{MainAlgorithm}, and then divide it into $\mc{N}=2,3,\ldots,$ $\mc{N}_m$ blocks using Algorithm \hyperlink{algo:2}{2}. If one wants to directly promote block-wise sparsity in $\bbK^*$ for any number of blocks $\mc{N}$, then one needs to know the corresponding initial block-structure a priori. We can use the block-structure obtained from Algorithm \hyperlink{algo:2}{2} for each $\mc{N}$ as that initial given structure and use the algorithm proposed in our recent paper \cite{negi2020sparsity} to carry out direct block-sparsity promotion in $\bbK^*$. The corresponding $\sbw$, $\scn$ and $J$ values obtained are compared with the output of Algorithm \hyperlink{algo:2}{2} in Fig. \ref{fig:caseb2BlockSparsePlot}. While $\sbw$ is held constant in Algorithm \hyperlink{algo:2}{2} (since it is already optimized in Algorithm \ref{MainAlgorithm}), the $\sbw$ for block-sparse algorithm increases steeply as $\mc{N}$ increases. On the other hand, $\scn$ for a given $\mc{N}$ is the same for both the algorithms because of the same corresponding block-structures. However, the $\mc{H}_2$ norm obtained through Algorithm \hyperlink{algo:2}{2} in each case is lower than that obtained through direct block-sparsity promotion for all $\mc{N}$. The reason is that the latter must collectively sparsify an entire block as a result of which the $\mc{H}_2$ norm becomes more conservative. Our element-wise sparsification approach in Algorithm \ref{MainAlgorithm} manages to avoid this conservatism.

\section{Conclusion}
\label{sec:conclusion}
This paper presented the co-design of network delays and sparse controller for LTI systems to improve their $\mc{H}_2$ performance. Bandwidth cost constraint is imposed to ensure finite bandwidth distribution among the communication links. The challenges of co-design borne out of implicit functional relationships between the delays, the sparse controller, and the $\mc{H}_2$-norm are overcome by a hierarchical algorithm, where the inner loop and outer loop are based on ADMM and SDP relaxations, respectively. Additional algorithms are derived for carrying out structural modifications in the sparse controller to manipulate the number and computation overhead of CNs such that both the $\mc{H}_2$ norm and the price of renting the agents are further reduced, at the cost of foregoing privacy of information. Numerical simulations show the effectiveness of the designs and bring out interesting observations about the relationship between the delays, sparsity, and $\mc{H}_2$ performance of the closed-loop system.

\section*{Appendix}

\subsection{Proof of Lemma \ref{lemmatau}}
The $ij$-th block of $\tilde{\bbA}$ is given as:
\begin{equation}
\tilde{\bbA}_{ij}=\begin{cases}
\sum\limits_{k=1, \ k \neq j}^N \frac{1}{\theta_j -\theta_k} \bb{I}_n, \hspace{1.7cm} i\in\mathbb{N}_{N-1}, i=j \\
\frac{1}{\theta_j - \theta_i}\prod\limits_{m=1, \ m\neq j, i}^N \frac{\theta_i - \theta_m}{\theta_j - \theta_m} \bb{I}_n, \ \ \ i\in\mathbb{N}_{N-1}, \ i\neq j\\
\bbA, \hspace{3.9cm} i=N, j\in\mathbb{N}_{N-1}. 
\end{cases} \label{Aorig}
\end{equation}
Substituting \eqref{theta} in \eqref{Aorig}, the diagonal and off-diagonal block matrices of the first $N-1$ block rows are given by:
\begin{subequations} 
\label{lambdaAvalue}
\begin{align}
&\tilde{\bbA}_{ii} =\frac{1}{\tau_o} \sum\limits_{k=1,\ k \neq i}^N  a_{ik}\bb{I}_n, \ \tilde{\bbA}_{ji} = \frac{a_{ji}}{\tau_o}\prod\limits_{m=1, \ m \neq j,i}^N \frac{a_{jm}}{a_{im}}\bb{I}_n,\\
&a_{ik}=-\left({ \sin\left(\frac{(2N-i-k)\pi)}{2(N-1)} \right) \sin \left(\frac{(k-i)\pi}{2(N-1)} \right)}\right)^{-1},
\end{align}
\end{subequations}
where $i,j\in \{1,\ldots,N\}$. Therefore, $\bs{\Lambda}$ can be written as:
\begin{align}
\bs{\Lambda}_{ij} =\begin{cases}
 \tilde{\bbA}_{ij}, \ \ \ \  \ i=1,\ldots,N-1, \ j=1,\ldots,N, \\
 \bb{0} , \ \ \ \ \ \ i=N, \ j=1,\ldots,N.
 \end{cases} \label{LambdaValue}
\end{align}
The proof follows from \eqref{Aorig}, \eqref{lambdaAvalue} and \eqref{LambdaValue}. \hfill $\blacksquare$
\smallskip

\subsection{Proof of Lemma \ref{lemmand}}
{From \eqref{spectral} and \eqref{acl}, $\bb{N}_d$ can be written as:
\begin{align}
\bb{N}_d = [l_1(-\tau_d),\ldots,l_N(-\tau_d)]^T \otimes \bb{I}_n. \label{Nd}
\end{align}
 Let $\vartheta_k = \cos \left(\frac{(N-k-1)\pi}{N-1}\right)$ for $k=\{0,\ldots,N-1\}$}. Using \eqref{theta}, \eqref{lj} and $c=\nicefrac{\tau_d}{\tau_o}$ we can write
\begin{align}
l_j(-\tau_d) = \prod\limits_{m=1, \ m\neq j}^N \frac{-c-0.5(\vartheta_{m-1} -1)}{0.5(\vartheta_{j-1} - \vartheta_{m-1})} .\label{prod}
\end{align} 
Using \eqref{Nd} and \eqref{prod}, $\bb{N}_d$ can be rewritten in the form of \eqref{Nddeclare} where the $j$-th row of $\bs{\Gamma}$ contains the coefficients of $l_j(-\tau_d)$. From \eqref{prod}, $l_j(-\tau_d)$ is a product of $N-1$ affine terms in $c$ whose coefficients are only dependent on $N$, and therefore, $\bs{\Gamma}$ is a constant for constant $N$. $\blacksquare$


\subsection{Proof of Theorem \ref{Theoremgrad}}
The proof of uniqueness of solution of \eqref{CC} and differentiability of $\bb{P}$ utilizes Lemma \ref{lemmatau} and \ref{lemmand}, and follows procedure similar to Theorem 2.1 and Lemma 3.1 in \cite{computational}, respectively. 
Specifically, $\bb{P}^{'}(\tau_o) \partial \tau_o$, $\bb{P}^{'}(c) \partial c$ and $\bb{P}^{'}(\bbK) dK$ follow as solutions of the following Lyapunov equations:
\begin{align}
&\bbA_{cl}^T \ \bb{P}^{'}(\tau_o) \partial \tau_o +\bb{P}^{'}(\tau_o) \partial \tau_o \ \bbA_{cl} = \frac{\partial \tau_o}{\tau^2_o} (\bs{\Lambda}^T\bb{P}+\bb{P}\bs{\Lambda} ) ,\label{taudiff}\\
& \bbA_{cl}^T \ \bb{P}^{'}(c) \partial c +\bb{P}^{'}(c) \partial c \ \bbA_{cl} = \bb{N}{'}_d \partial c\bb{K}^T_d\bb{G} +\bb{G}^T  \bbK_d \partial c^T {\bb{N}^{T}_d}^{'} , \label{cdiff}\\
&\bbA_{cl}^T\bb{P}^{'}(\bbK) \partial \bbK +\bb{P}^{'}(\bbK) \partial \bbK \ \bbA_{cl} =  -\bb{Z}_d - \bb{Z}^T_d - \bb{Z}_o - \bb{Z}^T_o,  \label{psiK}
 \end{align}
where $\bb{Z}_d=\bb{N}^T_d(\partial \bbK \circ\bs{\mc{I}}_d)\bb{G}$ and $\bb{Z}_o$ $=$ $\bb{N}^T_o$ $(\partial \bbK$ $\circ$ $\bs{\mc{I}}_o) \bb{G}$. The partial derivative of $J(\bbK)$ is
 $J' (\bbK) \partial \bbK = \te{Tr} ( \bs{\mc{B}}^T\bb{P}'$ $(\bbK)\bs{\mc{B}}) =\te{Tr}(\nabla J(\bbK)^T \partial \bbK)$, where $\partial \bbK\in\mathbb{R}^{m\times n}$. Post-multiplying \eqref{psiK} with $\bb{L}$ and taking its trace, we obtain
\begin{align}
\te{Tr}(\partial \bbK^T\nabla J(\bbK)) = \te{Tr} \big( \partial \bbK_d^T \bb{G} \bb{L} \bb{N}_d + \partial \bbK_o^T \bb{G}\bb{L} \bb{N}_o \big),\label{eq1}
\end{align}
where $\partial \bbK_d = \partial \bbK \circ \bs{\mc{I}}_d$ and $\partial \bbK_o = \partial \bbK \circ \bs{\mc{I}}_o$. Using the property $\te{Tr}((\bb{X}\circ \bb{Y})^T \bb{Z}) = \te{Tr}( \bb{X}^T (\bb{Y}\circ \bb{Z}))$ \cite[Prob. 8.37]{schott2016matrix}, where $\bb{X},\bb{Y},\bb{Z}\in\mathbb{R}^{m\times n}$ in \eqref{eq1}, we get \eqref{JK}. Using \eqref{taudiff} and \eqref{cdiff}, and a similar procedure as above, we obtain \begin{align}
& J'(\tau_o) = -\frac{1}{\tau_o^2} \te{Tr} ( \bs{\Lambda}^T \bb{P}\bb{L} + \bb{L} \bb{P}\bs{\Lambda} ), \label{33}\\
& J'(c) = \te{Tr}(\bb{N}^{'}_d \bbK_d^T \bb{G}\bb{L} +\bb{L} \bb{G}^T \bbK_d {\bb{N}^{'}_d}^T ). \label{34}
\end{align}
We can subsequently obtain \eqref{JtauJc} from \eqref{33} and \eqref{34}. $\blacksquare$


\subsection{Proof of Theorem \ref{theoremktau}}
Using $\bs{\phi}_0$, $\bs{\phi}_1$, $\psi_0$, $\bb{A}_1$ and $\Delta \tilde{\bb{C}}$ as stated in the theorem, we define $\bs{\phi}$ and $\psi$ using Lemma \ref{lemmatau} as:
\begin{align}
& \bs{\phi} = \bs{\phi}_0 + \bs{\phi}_1 + \bs{\phi}_2, \ \bs{\phi}_2= \bb{A}_1^T \Delta \bb{P} + \Delta \bb{P} \bb{A}_1^T, \nn \\
& \bs{\psi} = \bs{\psi}_0 + \bs{\psi}_1 , \  \bs{\psi}_1= \Delta\tilde{\bb{C}}^T \bb{R} \Delta \tilde{\bs{C}}.
\end{align}
The equation $\bs{\phi}+\bs{\psi} =0$ is equivalent to \eqref{CC} for $(\bbK,\omega_o,c^*)$ with $\omega_o = \nicefrac{1}{\tau_o}$ and therefore, $(\bbK,\omega_o,c^*)$ is a stabilizing tuple if $\bs{\phi}+\bs{\psi} \preceq 0$ is satisfied. This inequality will be satisfied by $\bs{\phi}$ and $\bs{\psi}$ if they satisfy $\lambda_{max}(\bs{\phi}) + \lambda_{max}(\bs{\psi}) \preceq 0$ \cite[Theorem 4.3.1 (Weyl)]{matrixanalysis}. Therefore, the following inequality is a sufficient condition for stability:
\begin{align}
&\bs{\phi}_0 + \bs{\phi}_1+ \bs{\psi}_0 + \lambda_{max} (\bs{\phi}_2) + \lambda_{max}(\bs{\psi}_1) \preceq 0. \label{eqbefore}
\end{align}
Equation \eqref{eqbefore} can be equivalently written as:
\begin{align}
\bs{\phi}_0 + \bs{\phi}_1+ \bs{\psi}_0  +\alpha I \preceq 0, \ \alpha \geq   \lambda_{max} (\bs{\phi}_2) + \lambda_{max}(\bs{\psi}_1).\label{lambdamax} 
\end{align}
\par \noindent Following \cite[Theorem 4.3.50]{matrixanalysis} and \cite[Theorem 1.2]{goldberg1982numerical}, $|\lambda_{max}$ $(\bs{\phi}_2)| \leq 2 \|\bbA^T_1 \Delta \bb{P}\|$, $|\lambda_{max}(\bs{\psi}_1)| = \|\bb{R}^{\nicefrac{1}{2}}\Delta \tilde{\bb{C}}\|^2$. Therefore, \eqref{theoremktaueq} yields the necessary $\alpha$ for satisfying \eqref{lambdamax}. $\blacksquare$
\subsection{Proof of Theorem \ref{theoremkc}}
Let $\bs{\phi}=\sum_{i=0}^4 \bs{\phi}_i$ where $\bs{\phi}_2=\bbA_1\Delta\bb{L} + \Delta\bb{L} \bbA^T_1$, $\bs{\phi}_3$ $= \bbA_2\bb{L}^* + \bb{L}^* \bbA^T_2$, $\bs{\phi}_4 =  \bbA_2 \Delta \bb{L} + \Delta \bb{L} \bbA^T_2$ and $\bbA_2 = -\bs{\mc{B}}\Delta \bbK_d \Delta \bb{N}^T_d$. The equation $\bs{\phi}+\bs{\mc{B}}\bs{\mc{B}}^T=0$ is equivalent to \eqref{LL} for $(\bbK,\tau^*_o,c)$, and $\bs{\phi}+\bs{\mc{B}}\bs{\mc{B}}^T\preceq 0$ implies that $(\bbK,\tau^*_o,c)$ is a stabilizing tuple. The rest of the proof can be obtained through similar arguments as Theorem \ref{theoremktau}.

\subsection{Derivation of (36)}
\label{subsec:derivationofkbar}
 Let $\thickbar{\bb{k}}=\texttt{vec}({\bbK} )$. Using the property $\texttt{vec}(\bb{ABC}) = (\bb{C}^T \otimes \bbA)\bbB$, on \eqref{kmineqgrad}, we obtain the following:
\begin{align}
&(\bs{\mc{T}}_{dd}(\thickbar{\bb{k}}\circ \bb{v}_d) + \bs{\mc{T}}_{od}(\thickbar{\bb{k}}\circ \bb{v}_o) )\circ \bb{v}_d + (\bs{\mc{T}}_{do}(\thickbar{\bb{k}}\circ \bb{v}_d) \nn\\
&+ \bs{\mc{T}}_{oo}(\thickbar{\bb{k}}\circ \bb{v}_o) )\circ \bb{v}_o + \rho (\thickbar{\bb{k}}\circ \bb{v}_d) + \rho (\thickbar{\bb{k}}\circ \bb{v}_o) = \mu. \label{63} \end{align}
Since $\bb{v}_d$ and $\bb{v}_o$ are binary vectors, $(\bs{\mc{T}}_{dd}(\thickbar{\bb{k}} \circ \bb{v}_d) )\circ \bb{v}_d = \big((\bs{\mc{T}}_{dd} \circ \hat{\bb{V}}^T_d)\thickbar{\bb{k}}\big)\circ \bb{v}_d$. Furthermore, $ \big((\bs{\mc{T}}_{dd} \circ \hat{\bb{V}}^T_d)\thickbar{\bb{k}}\big)\circ \bb{v}_d = (\hat{\bb{V}}_d \circ \bs{\mc{T}}_{dd} \circ \hat{\bb{V}}^T_d) \thickbar{\bb{k}}$. Substituting this in \eqref{63},
\begin{align}
&(\hat{\bb{V}}_d \circ \bs{\mc{T}}_{dd} \circ \hat{\bb{V}}^T_d + \hat{\bb{V}}_d \circ \bs{\mc{T}}_{od} \circ \hat{\bb{V}}^T_o +\hat{\bb{V}}_o  \nn\\
&\circ \bs{\mc{T}}_{do} \circ \hat{\bb{V}}^T_d + \hat{\bb{V}}_o \circ \bs{\mc{T}}_{oo} \circ \hat{\bb{V}}^T_o   + \rho I_{n^2}) \thickbar{\bb{k}}  = \mu.
\end{align}
We get \eqref{tddactual} from above, thereby completing the proof. $\blacksquare$


\subsection{Proof of Proposition \ref{propositionkcktau}}

1) The delay ratio $c = \frac{\tau_d}{\tau_o}$ can only be theoretically perturbed between $0$ and $1$ as $\tau_d\leq \tau_o$. Due to $\tau_{dpr}\neq 0$ and $\tau_{cpr}\neq 0$, $c$ can only be perturbed in the open interval $(\frac{\tau_{dpr}}{\tau_o^*}, 1-\frac{\tau_{cpr}}{\tau_o^*})$ if $\tau_o=\tau_o^*$ is kept constant. The minimum value of this interval is obtained when $\tau_{dtr}=0$ and the maximum value is obtained when $\tau_{ctr}=0$. 
\par \noindent Let $c^*$ be perturbed to $c\in$ $\left[{\tau_{dpr}}/{\tau_o^*},1-({ \tau_{cpr}}/{\tau_o^*}) \right]$ resulting in a cost $\sbw(c)$. Then, $\delta \sbw(c)$ is given as:
\begin{align}
 =& \sbw(c) - \sbw^*\\
    =& \frac{2 m_{cp} n_{cp}^*}{c \tau_o^* - \tau_{dpr}} + \frac{m_{cc}n_{cc}^* }{(1-c) \tau_o^* - \tau_{cpr}} - \sbw^* \\
     = & \dfrac{\splitfrac{\underbrace{\Big(\sbw^* {\tau_o^*}^2 \Big)}_{p_1} c^2 + \underbrace{\Big(   \tau_o^* \big(-\sbw^* (\tau_o^* + ( \tau_{dpr} - \tau_{cpr}) )- 2 m_{cp} n_{cp}^* + m_{cc} 
 \ncc^* \big)\Big)}_{q_1}c}{ + \underbrace{\Big(\big(\sbw^*\tau_{dpr} + 2m_{cp} n_{cp}^* \big)\big(\tau_o^* - \tau_{cpr}\big) - m_{cc}\ncc^* \tau_{dpr}\Big)}_{r_1} }}{\underbrace{\Big(-{\tau_o^*}^2\Big)}_{p_2}c^2 + \underbrace{\Big( \tau_o^* \big(\tau_o^* + (  \tau_{dpr}  - \tau_{cpr})  \big) \Big)}_{q_2}c + \underbrace{\Big(-\tau_{dpr} (\tau_o^* - \tau_{cpr})\Big)}_{r_2}}.
\end{align}
Thus, the constraint $\delta \sbw(c)\leq 0$ can be written as:
\begin{align}
(p_1- p_2) c^2 + (q_1 - q_2 ) c + (r_1-r_2) \leq  0 \ \implies  \ \delta \sbw(c) \leq  0. 
\end{align}
Since $p_1-p_2 = (\sbw^*+1) {\tau_o^*}^2 >0$, $\delta \sbw(c)\leq 0$ is a convex constraint w.r.t $c$.
\smallskip
\par \noindent 2) Keeping $c=c^* $ constant, we can write:
\begin{align}
    &c^* = \frac{\tau_d^*}{\tau_o^*}= c = \frac{\tau_d}{\tau_o} =  \frac{\tau_{dtr} + \tau_{dpr}}{\tau_{dtr} + \tau_{dpr} + \tau_{ctr} + \tau_{cpr}}, \\
    \implies  \ & \tau_{ctr} = \frac{\thickbar{c}^*(\tau_{dtr} + \tau_{dpr}) }{c^*} - \tau_{cpr}, \\
    \implies \ & c^* \tau_{ctr} - \thickbar{c}^* \tau_{dtr}  = \thickbar{c}^* \tau_{dpr} - c^* \tau_{cpr}. \label{75}
\end{align}
Thus, $\tau_o^*= \tau^*_{ctr} + \tau^*_{dtr} + \tau_{cpr} + \tau_{cpr}$ is perturbed to $\tau_o = \tau_{ctr} + \tau_{dtr} + \tau_{cpr} + \tau_{cpr}$ such that \eqref{75} is true. The constraint $\delta \sbw(\tau_o)$ is written as:
\begin{align}
   = \ & \sbw(\tau_o) - \sbw^*\\
   = \ & \frac{2 m_{cp} n_{cp}^*}{c^* \tau_o - \tau_{dpr}} + \frac{m_{cc}n_{cc}^* }{\thickbar{c}^* \tau_o - \tau_{cpr}} - \sbw^*\\
   = \ & \dfrac{\splitfrac{\underbrace{\Big(-\sbw^* c^*\thickbar{c}^* \Big)}_{p_3} \tau^2_o + \underbrace{\Big(2m_{cp}n^*_{cp}\thickbar{c}^* + m_{cc} n^*_{cc} c^* + 
\sbw^*(c^*\tau_{cpr}+\thickbar{c}^*\tau_{dpr}) \Big)}_{q_3} \tau_o}{+ \underbrace{\Big(-\sbw^* \tau_{dpr}\tau_{cpr}- 2m_{cp}n^*_{cp}\tau_{cpr} - m_{cc}n^*_{cc}\tau_{dpr}\Big)}_{r_3}}}{\underbrace{\Big(c^*\thickbar{c}^*\Big)}_{p_4}\tau^2_o + \underbrace{\Big( - c^* \tau_{cpr}-\thickbar{c}^*  \tau_{dpr}\Big)}_{q_4}\tau_o + \underbrace{\Big( \tau_{dpr}\tau_{cpr}\Big)}_{r_4}}.
\end{align}
Thus, the constraint $\delta\sbw(\tau_o) \leq 0$ can be written as:
\begin{align}
    (p_3-p_4) \tau_o^2 + (q_3-q_4) \tau_o + (r_3-r_4) \leq 0 \ \implies  \ \delta\sbw(\tau_o)\leq 0.
\end{align}
Since $p_3-p_4 = -(\sbw^* +1)c^* \thickbar{c}^*<0$, $\delta \sbw(\tau_o)\leq 0$ is a concave constraint w.r.t $\tau_o$.

\subsection{Proof of Lemma \ref{lemma:taucpr1}}
There are a total of $\mc{N}_m = \max(m,n)$ CNs available in the cyber-layer of a given CPS. To implement any topology $\bbT|_{\mc{N}}$ for a given $\bbK^*$, one needs to choose a subset of $\mc{N}$ CNs out of these $\mc{N}_m$ CNs. Let these subsets be denoted as $V|_{\mc{N}}$, $\mc{N}=\{1,\cdots,\mc{N}_m\}$. For example, any topology for $\mc{N}=4$ CNs will be implemented over the set $V|_{4}$. Let $\ell_{i,j} = \ell_{j,i}$ be the physical length of the edge between CN $i$ and $j$. We can then define a superset of physical lengths of the intra-layer links that will be utilized by a topology $\bbT|_{\mc{N}}$ as:
\begin{equation}
    \bs{\mf{L}}|_{\mc{N}} := \{\ell_{i,j} (= \ell_{j,i}) \ : \ i,j \in V|_{\mc{N}}  \}.
\end{equation}
Since for a given $\bbT|_{\mc{N}}$, we assume the worst case intra-layer link propagation delay as given in Sec. \ref{subsec:effectondelays}, $\tau_{cpr}(\bbT|_{\mc{N}})$ is a function of $\max(\bs{\mf{L}}|_{\mc{N}})$ written as:
\begin{equation}
    \tau_{cpr}(\bbT|_{\mc{N}}) = \bar{\kappa} \max(\bs{\mf{L}}|_{\mc{N}}), \ \bar{\kappa} =1/{\texttt{c}},
\end{equation}
where $\texttt{c}$ is the speed of light in the intra-layer link medium \citep{ugtext}. Let $\ell_{p,q} = \max \left(\bigcup\limits_{\mc{N}} \bs{\mf{L}}|_{\mc{N}}\right)$ be the largest link length between any two CNs. Thus, we can infer
\begin{align}
    \tau_{cpr}(\bbT^{(in)}|_{\mc{N}_m} )= \bar{\kappa} \max(\bs{\mf{L}}|_{\mc{N}_m}) = \bar{\kappa} \ \max \left(\bigcup_{\mc{N}} \bs{\mf{L}}|_{\mc{N}}\right)  = \ell_{p,q}.
\end{align}
One can choose the subset of $\mc{N}_m-1$ CNs to implement a topology $\bbT|_{\mc{N}_m-1}$ as $V|_{\mc{N}_m-1} = V|_{\mc{N}_m} / p$, i.e., all CNs except CN $p$ (or $q$) that corresponds to the source/destination of the link with the largest length in $\bs{\mf{L}}|_{\mc{N}_m}$. Therefore, the corresponding link length superset becomes $\max(\bs{\mf{L}}|_{\mc{N}_m-1}) \leq \ell_{p,q} $. Clearly, for this choice of $\mc{N}_m-1$ CNs, $\tau_{cpr}(\bbT|_{\mc{N}_m-1})\leq\tau_{cpr}(\bbT^{(in)}|_{\mc{N}_m}) $. We can apply a similar logic for any $\mc{N}<\mc{N}_m$ and choose $V|_{\mc{2}}$, $V|_{\mc{3}}$, $\cdots$,$V|_{\mc{N}_m-1}$, $V|_{\mc{N}_m}$ such that
\begin{align}
&\max(\bs{\mf{L}}|_{2}) \leq \max(\bs{\mf{L}}|_{3}) \leq \cdots \leq \max(\bs{\mf{L}}|_{\mc{N}_m-1}) \leq \max(\bs{\mf{L}}|_{\mc{N}_m}), \\
\implies  \    &\tau_{cpr}(\bbT|_{2}) \leq \tau_{cpr}(\bbT|_{3}) \leq \cdots \leq \tau_{cpr}(\bbT|_{\mc{N}_m-1}) \leq \tau_{cpr}(\bbT|_{\mc{N}_m}) = \tau_{cpr}(\bbT^{(in)}|_{\mc{N}_m}).
\end{align}

\subsection{Proof of Proposition \ref{prop:sec6}}
We use the shorthand $\sum \bb{a}|b$ for $\sum_{i=1}^\mc{N} \bb{a}_i (\bbK^*$, $\bbT^*|_b)$ and $\thickbar{b} = b-1$ for any $\bb{a}\in\mathbb{R}^\mc{N}$ and $b\in\mathbb{R}$ for this proof. Using \eqref{condition1}, for some $\mc{N}< {\mc{N}_m}$ we can write

\begin{align}
   & \frac{{\mc{N}_m} \thickbar{\mc{N}}_m }{\mc{N}\thickbar{\mc{N}}} \leq \frac{\sum \noff|_{{\mc{N}_m}}}{\sum\noff|_{\mc{N}} } \\
\implies \ & \frac{{\mc{N}_m}\thickbar{\mc{N}}_m}{\mc{N}\thickbar{\mc{N}}} \leq  \ \frac{\sum \noff|_{{\mc{N}_m}}\sum \bs{\mf{n}}|_{{\mc{N}_m}}}{\sum \noff|_{\mc{N}}\sum \bs{\mf{n}}|_{\mc{N}} }, \\
 \implies \ &   \frac{{\mc{N}_m}\thickbar{\mc{N}}_m \ncc|_{\mc{N}}}{\mc{N}\thickbar{\mc{N}}}  \leq  \ \ncc|_{\mc{N}_m} + \left(\sum\limits_i {n_{\text{off}}}_i|_{{\mc{N}_m}}\right) \left(\sum\limits_{j\neq i} n_j|_{{\mc{N}_m}}\right), \label{lastone}
\end{align}
where $\ncc|_{\mc{N}} = {\bs{\mf{n}}|_{\mc{N}}}^T \noff|_{\mc{N}}$. From \eqref{lastone}, $\ncc|_{\mc{N}} \leq \ncc|_{\mc{N}_m}$ is true for some ${\mc{N}} <\mc{N}_m$ for any given $\noff|_{{\mc{N}_m}}$ and $\bs{\mf{n}}|_{{\mc{N}_m}}$ if
\begin{equation}
\frac{{\mc{N}_m}\thickbar{\mc{N}}_m}{\mc{N}\thickbar{\mc{N}}}  \geq  \ 1 + \frac{\left(\sum {n_{\text{off}}}_i|_{{\mc{N}_m}} \sum_{j\neq i} n_j|_{{\mc{N}_m}}\right)}{\ncc|_{\mc{N}_m}}. \label{maximize1}
\end{equation}
The maximum value of RHS in \eqref{maximize1} is $n$; we prove this in \textbf{(C3)} shortly. Therefore, we can rewrite \eqref{maximize1} as:
\begin{align}
   \frac{{\mc{N}_m}\thickbar{\mc{N}}_m}{\mc{N}\thickbar{\mc{N}}}  \geq  \ n. \label{neweq}
\end{align}
We will now prove the theorem separately for the cases $n\leq m$ and $n>m$.
\par \noindent \textbf{(C1) Case $n\leq m$} \par If the following constraints are satisfied for any $n,m$ with $n\leq m$:
\begin{align}
  n \geq  \ {\mc{N}}({\mc{N}}-1) + 1, \ \   {\mc{N}} \in  \ [2,n], \label{case1sol} 
  \end{align}
then, from \eqref{neweq}, the theorem is true. Clearly the above constraints are satisfied for any $n\geq 2$. Therefore, the theorem is proven for the case $n\leq m$.

\par \noindent \textbf{(C2) Case $n> m$} \par Since $n>m$, ${\mc{N}_m}=m$ and ${\mc{N}}\in [2,m]$. We can thus rewrite \eqref{maximize1} as:
\begin{align}
     \frac{m(m-1)}{{\mc{N}}({\mc{N}}-1)}  \geq &  \ n. \label{15}
\end{align}
Therefore, using \eqref{15} and simple algebra, we can conclude that the above constraints will be satisfied simultaneously, i.e., the theorem is valid only when $m \leq  n \leq ({m(m-1)})/{2}$ for the case $n>m$.
\par \noindent \textbf{(C3) Proof for the maximum value of the RHS in \eqref{maximize1}}
\par We need to prove that the maximum value of RHS of \eqref{maximize1} is $n$. To do that, we write the RHS of \eqref{maximize1} as $1+\phi$ and prove that the maximum value of $\bs{\phi}=n-1$, i.e.,
\begin{equation}
    \underset{{\mc{N}_m}\geq 2}{\text{max}} \ \ \phi= \left(\frac{\sum_i \noffs_i|_{{\mc{N}_m}} \sum_{j, \ j\neq i} n_j|_{{\mc{N}_m}}}{ {\bs{\mf{n}}|_{{\mc{N}_m}}^T}{\noff|_{{\mc{N}_m}}}}\right) = n-1, \label{87}
\end{equation}
given that $\sum_i n_i|_{{\mc{N}_m}} = n$, $n_i|_{{\mc{N}_m}} \in [1,n-\mc{N}_m]$ and $ \noffs_i|_{{\mc{N}_m}}\in [0,{\mc{N}_m}-1]$.  For ease of notation, we write $\noff|_{\mc{N}_m} = [\noffs_1,\noffs_2,\ldots,\noffs_{\mc{N}_m}]^T$ and $\bs{\mf{n}}|_{\mc{N}_m} = [n_1,n_2,\ldots,n_{\mc{N}_m}]^T$. Substituting $n_{\mc{N}_m} =n - \sum_i n_i$ in \eqref{87}, we can write
\begin{align}
 \phi = \frac{\big(\sum_i^{\mc{N}_m-1} \noffs_i (n-n_i) \big) + \noffs_{\mc{N}_m} \sum_i^{\mc{N}_m-1} n_i }{\sum_{i=1}^{\mc{N}_m-1} \noffs_i n_i + n \noffs_{\mc{N}_m} - \noffs_{\mc{N}_m} \sum_{i=1}^{\mc{N}_m-1} n_i }.
\end{align}
Let $\Delta \noffs_i = \noffs_{\mc{N}_m} - \noffs_i$ $\forall$ $i$. Then, we can write:
\begin{equation}
    \phi = \frac{(\sum_{i=1}^{\mc{N}_m-1} \noffs_i )n + \sum_{i=1}^{\mc{N}_m-1} n_i\Delta \noffs_i  }{ n \noffs_{\mc{N}_m} - \sum_{i=1}^{\mc{N}_m-1} n_i \Delta \noffs_i}. \label{89}
\end{equation}
Substituting $\noffs_i = \noffs_{\mc{N}_m} - \noffs_i$ $\forall$ $i\in\mathbb{N}_{\mc{N}_m}$ in \eqref{89}:
\begin{equation}
 \phi = \frac{(\mc{N}_m -1) n \noffs_{\mc{N}_m} + \sum_{i=1}^{\mc{N}_m-1} (n_i-n)\Delta \noffs_i  }{n \noffs_{\mc{N}_m} - \sum_{i=1}^{\mc{N}_m-1} (n_i-n)\Delta \noffs_i  }. \label{90}
\end{equation}
To maximize ${\phi}$, we need to choose $\noffs_i$, $i\in\mathbb{N}_{\mc{N}_m}$ such that we maximize $\sum_{i=1}^{\mc{N}_m-1}$ $(n_i-n)$ $\Delta \noffs_i $ and $\sum_{i=1}^{\mc{N}_m-1}$ $n_i$  $\Delta \noffs_i$, which is only possible if we choose $n_i$ and $\Delta \noffs_i$ $\forall$ $i\in\mathbb{N}_{\mc{N}_m-1}$ to be maximum. Therefore, we assume in \eqref{90} that $\sum_{i=1}^{\mc{N}_m-1} n_i = n-1$ and $\Delta \noffs_i = \mc{N}_m-1$ $\forall$ $i\in\mathbb{N}_{\mc{N}_m-1}$. Substituting these in \eqref{90}:
\begin{align}
    &\max {\phi} = \frac{n(\mc{N}_m-1)^2  + (\mc{N}_m-1)( (\sum_{i=1}^{\mc{N}_m-1}n_i) - (\mc{N}_m-1)n) }{n(\mc{N}_m-1) - (\sum_{i=1}^{\mc{N}_m-1}n_i)  (\mc{N}_m-1) }\\
&   \max  \ {\phi} = \frac{(\mc{N}_m-1)n  + (n-1 - (\mc{N}_m-1)n) }{n-(n-1)} = n-1. \nn
\end{align}
Therefore, \eqref{87} is proved. Since RHS of \eqref{maximize1} is $1+\phi$, its maximum value is $n$.
\begin{figure}[ht]
    \centering
    \includegraphics[scale=0.5]{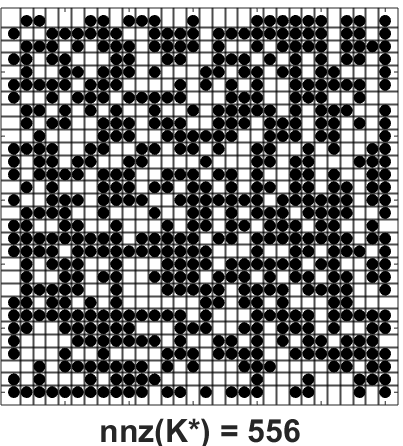}
    \caption{Block Structure $\bmcK(\bbK^*,\bbT^{(in)}|_{\mc{N}_m=30})$}
    \label{fig:Kinitial}
\end{figure}
\subsection{Initial Conditions for Example in Section \ref{sec:simulations2}}
\label{subsec:initialconditions}
The obtained sparse $\bbK^*$ from Algorithm \ref{MainAlgorithm} for Case B has $\nnz(\bbK^*)=556$. The initial topology $\bbT^{(in)}|_{\mc{N}_m=30}$ is randomly chosen as:
\begin{align}
&\bs{\mf{X}}^{(in)}|_{\mc{N}_m} = \{     30, \    20, \    22, \     1, \    25, \    18, \    11, \    24, \    16, \     2, \    28, \    26, \     3, \     5, \     7, \     4, \    10, \    12, \     6, \    21, \    27, \     9, \    15, \    19,  \nn \\ 
 & \hspace{12cm}     8, \ 17, \    23, \    13, \    29, \    14 \} \\
  &\bs{\mf{U}}^{(in)}|_{\mc{N}_m} = \{  26, \   15, \    24, \    12, \     9, \     1, \     13, \    27, \     8, \    10, \    22, \    11, \     5, \    23, \    16, \    20, \     6, \    14, \    19, \    25, \     7, \    17, \     4, \    21,\nn\\
   & \hspace{11cm}     18, \     2, \     3, \   29, \    28, \    30 \} \\
   &\bs{\mf{n}}^{(in)}|_{\mc{N}_m} = [1,1,1,\cdots, 1,1]^T \in \mathbb{Z}^{30}, \ \  \bs{\mf{m}}^{(in)}|_{\mc{N}_m} = [1,1,1,\cdots, 1,1]^T \in \mathbb{Z}^{30}
\end{align}
The corresponding block structure $\bmcK(\bbK^*,\bbT^{(in)}|_{\mc{N}_m=30})$ is given in Fig. \ref{fig:Kinitial}.

\end{document}